\documentclass[11pt,a4paper,reqno]{article}
\usepackage{mysty}
\usepackage{authblk}

\usepackage{fullpage}
\usepackage{subfig}

\usepackage{amsfonts,amssymb}
\usepackage{amsmath}
\usepackage{graphicx}
\usepackage{cite}
\usepackage{enumerate}

\theoremstyle{plain}

\usepackage{babel}
\providecommand{\theoremname}{Theorem}

\SetLabelAlign{parright}{\parbox[t]{\labelwidth}{\raggedleft{#1}}}
\setlist[description]{style=multiline,topsep=4pt,align=parright}

\makeatletter
\let\reftagform@=\tagform@
\def\tagform@#1{\maketag@@@{(\ignorespaces\textcolor{black}{#1}\unskip\@@italiccorr)}}
\newcommand{\iref}[1]{\textup{\reftagform@{\tcr{\ref{#1}}}}}
\makeatother

\begin{document}
	
\title{A Kernel-free Boundary Integral Method for the Bidomain Equations}
\author{Xindan Gao, Li Cai, Craig S. Henriquez and Wenjun Ying}

\affil{School of Mathematical Sciences and Institute of Natural Sciences, Shanghai Jiao Tong University, Shanghai, China}

\date{}
\maketitle

\begin{abstract}
	The bidomain equations have been widely used to mathematically
	model the electrical activity of the cardiac tissue. In this 
	work, we present a potential theory-based Cartesian grid method 
	which is referred as the kernel-free boundary integral (KFBI) 
	method which works well on complex domains to efficiently simulate the linear diffusion part of 
	the bidomain equation. After a proper temporal discretization, 
	the KFBI method is applied to solve the resulting homogeneous 
	Neumann boundary value problems with a second-order accuracy. 
	According to the potential theory, the boundary integral 
	equations reformulated from the boundary value problems can 
	be solved iteratively with the simple Richardson iteration or
	the Krylov subspace iteration method. During the iteration, 
	the boundary and volume integrals are evaluated by limiting 
	the structured grid-based discrete solutions of the 
	equivalent interface problems at quasi-uniform interface nodes
	without the need to know the analytical expression of Green's 
	functions. In particular, the discrete linear system of the 
	equivalent interface problem obtained from the standard 
	finite difference schemes or the finite element schemes
	can be efficiently solved by fast elliptic solvers such 
	as the fast Fourier transform based solvers or those based 
	on geometric multigrid iterations after an appropriate
	modification at the irregular grid nodes. Numerical results 
	for solving the FitzHugh-Nagumo bidomain equations in both 
	two- and three-dimensional spaces are presented to
	demonstrate the numerical performance of the KFBI method such as 
	the second-order accuracy and the propagation and scroll wave
	of the voltage simulated on the real human left ventricle model.
	
\end{abstract}

\section{Introduction}

The set of bidomain equations, consisting of the equations for
the intra- and extracellular potentials, is currently the most 
complete mathematical model for describing the electrical 
activities of the cardiac tissue \cite{clayton2008guide}. 
The cardiac simulations are inherently computationally expensive 
on account of several factors such as the high resolution in 
space and time discretization and the high complexity of living 
systems. To be specific, the numerical computational time is 
mainly consumed in a burdensome part resulting from 
the linear system discretized from the bidomain equations
which should be solved repeatedly \cite{vigmond2008solvers}. 

In detailed numerical simulations, for the temporal discretization, 
the explicit methods such as the forward Euler method and Runge 
Kutta method are the easiest way to implement. Although they do 
not require matrix inversion, the timesteps allowed in these methods are 
generally very small to satisfy the requirement of stability 
\cite{dos2004parallel,potse2006comparison,quan1998efficient}. 
The fully implicit methods, preferred for stiff systems, have 
no timestep restriction, but they often have the significant 
drawback of requiring the numerical solution of a very 
large-scale nonlinear system at each timestep
\cite{hooke1993efficient,murillo2004fully,ying2008efficient}.
The semi-implicit method combines the advantages of the 
above methods that it allows for a much larger timestep than 
that used in the explicit method and requires the numerical solution 
of a linear system at each timestep \cite{keener1998numerical}.
Operator splitting such as the first-order Godunov splitting 
and the second-order Strang splitting is a popular technique being 
able to uncouple the cell-model ordinary differential equation 
(ODE) system from the electro-diffusive parts of the bidomain equations
\cite{keener1998numerical,qu1999advanced,whiteley2006efficient,sundnes2001efficient}.
In this work, the second-order Strang splitting is 
applied to split the bidomain equations into two 
manageable parts, $i.e.$, the nonlinear reaction
and linear diffusion parts.

For the spatial discretization of the domain, the finite 
difference method (FDM) \cite{clayton2008guide,pollard1992cardiac,trew2005generalized}
is usually applied to solve the bidomain equations on a 
regular computational space, but it lacks the ability   
to follow complex surfaces. Although the finite element 
method (FEM) \cite{fischer2000bidomain,sundnes2007computing,trayanova2006defibrillation}
and the finite volume method (FVM) \cite{jacquemet2005finite,trew2005finite}
have the advantages of being able to model boundary conditions 
at the curved surfaces of the heart, they are generally more 
difficult to implement and more computationally expensive than the FDM. 

In our implementation, after applying the second-order 
Strang splitting to decouple the bidomain equations 
into the nonlinear reaction and linear diffusion parts, 
we discrete the nonlinear reaction part by the forward 
Euler method in the first half step and by the backward 
Euler method in the second half step. The resulting discrete 
nonlinear systems are space independent which can be 
solved efficiently. We integrate the linear diffusion part
by the second-order mid-point method where the obtained 
semi-discrete equation is a coupled elliptic system 
with a possible anisotropy. In this work, we present an 
accurate and efficient algorithm, a generalized boundary integral method,
to solve the semi-discrete diffusion part of the bidomain 
equations, $i.e.$, the so-called kernel-free boundary integral 
(KFBI) method. As a direct extension of the grid-based 
boundary method by Mayo \cite{mayo1984fast,mayo1985fast,mayo1992rapid}, 
the KFBI method proposed by Ying $et.al.$ \cite{ying2007kernel,ying2014kernel}
can efficiently solve the variable coefficients 
elliptic partial differential equations (PDE) with possible
anisotropy and inhomogeneity.

When solving the coupled semi-discrete diffusion system with a Neumann 
boundary condition by the KFBI method with second-order accuracy, we 
first reformulate the system as the boundary integral equations (BIEs) 
where the density defined on the domain boundary 
is unknown. The BIEs then can be iteratively solved 
with a Krylov subspace iteration method such as
the generalized minimal residual (GMRES) \cite{saad1986gmres} 
iteration or the simple Richardson iteration.
During the iteration, the boundary and volume integrals 
can be evaluated by computing the limit values of the structured grid-based 
numerical solutions of the equivalent interface 
problem without the need to know the analytical expressions 
of Green's functions. Thus, the KFBI method is said to be kernel-free.
After solving an equivalent interface problem on a Cartesian grid by 
the FDM or the FEM, we quadratically interpolate the discrete solution 
at the discretization nodes of the domain boundary to evaluate the 
values of the boundary and volume integrals. We point out that
the KFBI method can preserve the symmetric and positive definite 
property of the coefficient matrix of the discrete linear system 
of the equivalent interface problem. Hence, the 
discrete linear system can be efficiently solved with the
standard geometric multigrid iterations or the fast Fourier
transform based solvers. To show the numerical performance of the 
KFBI method, we simulate the bidomain equations on different complex
regions such as the heart-shaped domain and a real 
human left ventricle (LV) space model in two- and 
three-dimensional spaces, respectively.

\section{The Model}

The bidomain equation which is governed by a singularly perturbed 
reaction-diffusion system consists of a set of nonlinear ordinary 
differential equations representing cell membrane dynamics and the 
partial differential equations representing the propagation of the 
electrical signal through the cardiac tissue \cite{keener1998mathematical}.
Let $\Omega\subset\mathbb{R}^{d}$ $(d=2 \textrm{ or } 3)$ be the 
bounded region occupied by the cardiac tissue which is usually 
complicated. We begin with the bidomain equation which consists
of the equations for the intra- and extracellular potentials, 
$\Phi_{\textrm{i}}$ and $\Phi_{\textrm{e}}$, coupled through 
the transmembrane potential $V_{\textrm{m}}=\Phi_{\textrm{i}}-\Phi_{\textrm{e}}$,
for $\mathbf{x}\in\Omega$ and $t>0$

\begin{equation}
\begin{aligned}C_{\text{m}}\dfrac{\partial V_{\textrm{m}}}{\partial t}+I_{\textrm{ion}}(V_{\textrm{m}},\mathbf{q}) & =\dfrac{1}{\beta}\nabla\cdot(\mathbf{D}_{\textrm{i}}\nabla\Phi_{\textrm{i}}),\\
C_{\textrm{m}}\dfrac{\partial V_{\textrm{m}}}{\partial t}+I_{\textrm{ion}}(V_{\textrm{m}},\mathbf{q}) & =-\dfrac{1}{\beta}\nabla\cdot(\mathbf{D}_{\textrm{e}}\nabla\Phi_{\textrm{e}})-I_{\text{stim}},\\
\dfrac{\partial\mathbf{q}}{\partial t} & =\mathcal{M}(V_{\textrm{m}},\mathbf{q}),
\end{aligned}
\end{equation}
where $C_{\textrm{m}}$ is the membrane capacitance per unit area; 
$\beta$ is a surface to volume ratio of the cardiac cells; 
$\mathbf{D}_{\textrm{i}}$ and $\mathbf{D}_{\textrm{e}}$ are the 
space dependent intracellular and extracellular conductivity tensors, 
respectively; $\mathbf{q}$ is a set of state variables; 
$I_{\textrm{ion}}(V_{\textrm{m}},\mathbf{q})$ and $\mathcal{M}(V_{\textrm{m}},\mathbf{q})$ 
are two known functions approximating the cellular membrane dynamics;
$I_{\textrm{stim}}=I_{\textrm{stim}}(t,\mathbf{x})$ is a given 
extracellular stimulus current. In this work, we consider the 
case that the cardiac tissue is insulated, $i.e.$, the bidomain 
equation subject to the homogeneous Neumann boundary condition

\begin{equation}
\begin{array}{ccccc}
\mathbf{n}\cdot(\mathbf{D}_{\textrm{i}}\nabla\Phi_{\textrm{i}})=0 & \textrm{and} & \mathbf{n}\cdot(\mathbf{D}_{\textrm{e}}\nabla\Phi_{\textrm{e}})=0 & \textrm{on} & \partial\Omega,\end{array}\label{eq:bdry cond}
\end{equation}
where $\mathbf{n}$ is the unit outward normal to the boundary 
$\partial\Omega$ and imposed on $\partial\Omega$ all the time. 
Provided some appropriate initial conditions on the intracellular 
potential $\Phi_{\text{i}}$, the extracellular potential $\Phi_{\text{e}}$, 
and the vector of the state variables $\mathbf{q}$, the bidomain 
equations can be numerically solved with some specific spatial discretization
and temporal integration methods.

\section{Operator Splitting Techniques}

Operator splitting techniques \cite{keener1998numerical,qu1999advanced,trangenstein2004operator}
are widely used to efficiently solve the bidomain equations by 
decoupling the reaction and diffusion parts in the reaction-diffusion 
system to be solved. The application of an operator splitting technique 
allows integrating the reaction and diffusion parts independently and 
implicitly without solving a large nonlinear system in each time step. 
Let $\Delta t>0$ be the timestep and $t^{n}=n\Delta t$ be the discrete 
times. The second-order Strang splitting is applied to numerically 
simulate the electrical activity of the cardiac tissue by integrating 
the decoupled parts in a symmetric way in the following three steps.

step 1: First half time step nonlinear reaction

\begin{equation}
\begin{aligned}C_{\text{m}}\dfrac{\partial V_{\text{m}}}{\partial t} & =-I_{\text{ion}}(V_{\text{m}},\mathbf{q}) & \textrm{for}\;t^{n}<t<t^{n}+\dfrac{\Delta t}{2},\\
\dfrac{\partial\mathbf{q}}{\partial t} & =\mathcal{M}(V_{\text{m}},\mathbf{q}) & \textrm{for}\;t^{n}<t<t^{n}+\dfrac{\Delta t}{2}.
\end{aligned}
\label{eq:ode1}
\end{equation}

step 2: A full time step linear diffusion

\begin{equation}
\begin{aligned}C_{\text{m}}\dfrac{\partial V_{\text{m}}}{\partial t} & =\dfrac{1}{\beta}\nabla\cdot(D_{\text{i}}\nabla\Phi_{\text{i}}) & \textrm{for}\;t^{n}<t<t^{n+1},\\
C_{\text{m}}\dfrac{\partial V_{\text{m}}}{\partial t} & =-\dfrac{1}{\beta}\nabla\cdot(D_{\text{e}}\nabla\Phi_{\text{e}})-I_{\text{stim}} & \textrm{for}\;t^{n}<t<t^{n+1}.
\end{aligned}
\label{eq:diffusion}
\end{equation}

step 3: Second half time step nonlinear reaction

\begin{equation}
\begin{aligned}C_{\text{m}}\dfrac{\partial V_{\text{m}}}{\partial t} & =-I_{\text{ion}}(V_{\text{m}},\mathbf{q}) & \textrm{for}\;t^{n}+\dfrac{\Delta t}{2}<t<t^{n+1},\\
\dfrac{\partial\mathbf{q}}{\partial t} & =\mathcal{M}(V_{\text{m}},\mathbf{q}) & \textrm{for}\;t^{n}+\dfrac{\Delta t}{2}<t<t^{n+1}.
\end{aligned}
\label{eq:ode2}
\end{equation}

Here, Eqs. (\ref{eq:ode1}) and (\ref{eq:ode2}) represent the resulting
nonlinear reaction equations which are space independent and usually are 
stiff \cite{sundnes2002numerical,wanner1996solving} while Eq. (\ref{eq:diffusion}) 
represents the resulting linear diffusion equation which is space dependent. 
After the application of the Strang operator splitting, the reaction and 
diffusion equations can be separately discretized by implicit integration 
methods such as the backward Euler method, semi-implicit method, and the 
Crank-Nicolson to satisfy the stability requirement without solving a large 
nonlinear system \cite{sundnes2005operator}. To be specific, during each 
timestep, we integrate the reaction equations (\ref{eq:ode1}) and (\ref{eq:ode2}) 
by the forward Euler method and the backward Euler method, respectively,
and integrate the diffusion equation (\ref{eq:diffusion}) by the second-order 
implicit midpoint scheme.

For ease of illustration, we first give the implicit midpoint scheme 
for the abstract form
\begin{equation}
\dfrac{du(t)}{dt}=f(u),\textrm{ for } t>0,\label{eq:evolution eq}
\end{equation}
which can then be straightforwardly extended to the bidomain equations.
The implicit midpoint scheme for Eq. (\ref{eq:evolution eq}) is given by
\[
\dfrac{u^{n+1}-u^{n}}{\Delta t}=f\left(\dfrac{u^{n+1}+u^{n}}{2}\right),
\]
and can be rewritten as the discrete format obtained from the application 
of the backward Euler method

\[
\dfrac{w-u^{n}}{\Delta t/2}=f(w),
\]
where $w=(u^{n+1}+u^n)/2$ and the corresponding timestep is $\Delta t/2$. 
In our implementation, we first obtain the numerical solution of the 
intermediate variable $w$ and then compute the numerical solution  $u^{n+1}=2w-u^n$. Therefore, we
restrict ourselves to the discussion on the application of the backward
Euler method to integrate the diffusion equation (\ref{eq:diffusion}). 
Time integration of the diffusion equation by the backward Euler method 
from time $t^{n}$ to time $t^{n+1}$ leads to the semi-discrete partial 
differential equations

\begin{equation}
\begin{aligned}C_{\textrm{m}}\dfrac{V_{\text{m}}^{n+1}-V_{\text{m}}^{n}}{\Delta t} & =\dfrac{1}{\beta}\nabla\cdot(\mathbf{D}_{\text{i}}\nabla\Phi_{\text{i}}^{n+1}),\\
C_{\text{m}}\dfrac{V_{\text{m}}^{n+1}-V_{\text{m}}^{n}}{\Delta t} & =-\dfrac{1}{\beta}\nabla\cdot(\mathbf{D}_{\text{e}}\nabla\Phi_{\text{e}}^{n+1})-I_{\text{stim}},
\end{aligned}
\label{eq:pde0}
\end{equation}
where $V_{\text{m}}^{n}$, $\Phi_{\text{i}}^{n}$, and $\Phi_{\text{e}}^{n}$ are 
the finite difference approximations of the variables $V_{\text{m}}$, 
$\Phi_{\text{i}}$, and $\Phi_{\text{e}}$ at time $t^{n}$, respectively. 

Assume that the extracellular stimulus current is only space dependent. 
We rewrite the semi-discrete partial differential equations (\ref{eq:pde0}) as

\begin{equation}
\begin{aligned}\nabla\cdot(\mathbf{D}_{\text{i}}\nabla\Phi_{\text{i}}^{n+1})-\kappa(\Phi_{\text{i}}^{n+1}-\Phi_{\text{e}}^{n+1}) & =-\kappa V_{\text{m}}^{n},\\
\nabla\cdot(\mathbf{D}_{\text{e}}\nabla\Phi_{\text{e}}^{n+1})+\kappa(\Phi_{\text{i}}^{n+1}-\Phi_{\text{e}}^{n+1}) & =\kappa V_{\text{m}}^{n}-\beta I_{\text{stim}},
\end{aligned}
\label{eq:pde1}
\end{equation}
where $\kappa=C_{\text{m}}\beta/\Delta t$.  Denote the unknown 
vector function and the right hand side of the system (\ref{eq:pde1}) 
by $\mathbf{u}\equiv(\Phi_{\text{i}}^{n+1},\Phi_{\text{e}}^{n+1})^{T}$
and $\mathbf{f}\equiv(-\kappa V_{\text{m}}^{n},\kappa V_{\text{m}}^{n}-\beta I_{\text{stim}})^{T}$.
Define Eq. (\ref{eq:pde1}) in the form of modified Helmholtz
equation as followed

\begin{equation}
\begin{array}{ccc}
(\mathcal{L}\mathbf{u})\equiv\nabla\cdot(\mathbf{D}\nabla\mathbf{u})-\mathbf{K}\mathbf{u}=\mathbf{f} & \textrm{in} & \Omega,\end{array}\label{eq:bvp}
\end{equation}
where $\mathbf{D}\equiv\left(\begin{array}{cc}
\mathbf{D}_{\text{i}} & 0\\
0 & \mathbf{D}_{\text{e}}
\end{array}\right)$ and $\mathbf{K}\equiv\left(\begin{array}{cc}
\kappa & -\kappa\\
-\kappa & \kappa
\end{array}\right)$.
The corresponding boundary condition is homogeneous Neumann boundary 
condition which is noted as 

\begin{equation}
\mathbf{g}^{N}\equiv\mathbf{n}\cdot(\mathbf{D}\nabla\mathbf{u})\equiv\left(\mathbf{n}\cdot(\mathbf{D}_{\text{i}}\nabla\Phi_{\text{i}}),\mathbf{n}\cdot(\mathbf{D}_{\text{e}}\nabla\Phi_{\text{e}})\right)^{T}=\mathbf{0}.\label{eq:BC}
\end{equation}
In the following sections \ref{sec:BIE} and \ref{sec:KFBI}, we 
introduce how to solve the boundary value problem (\ref{eq:bvp}) 
by the KFBI method during each timestep. 

\section{Boundary Integral Equation Formulation\label{sec:BIE}}

In this section, we reformulate the homogeneous Neumann boundary
value problem (\ref{eq:bvp}) as a BIE \cite{ying2007kernel}.
Let $\mathcal{B}$ be a larger regular domain which completely contains
the computational domain $\Omega$ and $\partial\Omega$ be the smooth
interface in $\mathcal{B}$ which separates the larger regular domain $\mathcal{B}$
into two subdomains $\Omega$ and $\Omega^{C}$. Let $G(\mathbf{q};\mathbf{p})$ 
be Green's function associated with the operator defined in Eq. (\ref{eq:bvp}) 
on the larger regular domain $\mathcal{B}$, which satisfies

\begin{equation}
\begin{aligned}\mathcal{L}G(\mathbf{q};\mathbf{p})\equiv\nabla_{\mathbf{q}}\cdot(\mathbf{D}\nabla_{\mathbf{q}}G)-\mathbf{K}G=\delta(\mathbf{q}-\mathbf{p}) &  & \mathbf{q}\in\mathcal{B},\\
\mathbf{n}_{\mathbf{q}}\cdot(\mathbf{D}\nabla_{\mathbf{q}}G(\mathbf{q};\mathbf{p}))=0 &  & \mathbf{q}\in\partial\mathcal{B},
\end{aligned}
\end{equation}
for each fixed $\mathbf{p}\in\mathcal{B}$. Here, $\delta(\mathbf{q}-\mathbf{p})$
is the Dirac delta function; $\nabla_{\mathbf{q}}$ stands for the
gradient operator with respect to the space variable $\mathbf{q}\in\mathbb{R}^{d}$.
In terms of the Green's function $G(\mathbf{q};\mathbf{p})$, the
solution to the Neumann boundary value problem defined by Eqs. (\ref{eq:bvp})
and (\ref{eq:BC}) can be expressed as a sum of a volume integral
and a boundary integral

\begin{equation}
\begin{array}{ccc}
\mathbf{u}(\mathbf{p})=\int_{\Omega}G(\mathbf{q};\mathbf{p})\mathbf{f}(\mathbf{q})d\mathbf{q}-\int_{\partial\Omega}G(\mathbf{q};\mathbf{p})\mathbf{\mathbf{\boldsymbol{\psi}}}(\mathbf{q})ds_{\mathbf{q}}, &  & \mathbf{p}\in\Omega,\end{array}\label{eq:solution}
\end{equation}
where the density $\boldsymbol{\psi}(\mathbf{p})$ satisfies the second kind Fredholm BIE

\begin{equation}
\begin{array}{ccc}
\dfrac{1}{2}\boldsymbol{\psi}(\mathbf{p})-\int_{\partial\Omega}\mathbf{n}_{\mathbf{p}}\cdot\mathbf{D}\nabla_{\mathbf{p}}G(\mathbf{q};\mathbf{p})\boldsymbol{\psi}(\mathbf{q})ds_{\mathbf{q}}=\mathbf{g}^{N}-\mathbf{n_{p}}\cdot\mathbf{D}\nabla_{\mathbf{p}}\left(\int_{\Omega}G(\mathbf{q};\mathbf{p})\mathbf{f}(\mathbf{q})d\mathbf{q}\right), &  & \mathbf{p}\in\partial\Omega.\end{array}\label{eq:bie}
\end{equation}

We can solve the BIE (\ref{eq:bie}) by the following simple iterations

\begin{equation}
\begin{aligned}
\text{\ensuremath{\boldsymbol{\psi}}}_{\nu+1}(\mathbf{p}) 
& =\boldsymbol{\psi}_{\nu}(\mathbf{p})+2\gamma\left[\mathbf{g}^{N}-\mathbf{n_{p}}\cdot\mathbf{D}\nabla_{\mathbf{p}}\left(\int_{\Omega}G(\mathbf{q};
\mathbf{p})\mathbf{f}(\mathbf{q})d\mathbf{q}\right)\right.\\ 
& \left.-\left(\dfrac{1}{2}\boldsymbol{\psi}_{\nu}(\mathbf{p})-\int_{\partial\Omega}\mathbf{n}_{\mathbf{p}}\cdot\mathbf{D}\nabla_{\mathbf{p}}G(\mathbf{q};\mathbf{p})\boldsymbol{\psi}_{\nu}(\mathbf{q})ds_{\mathbf{q}}\right)\right],\mathbf{p}\in\partial\Omega  
\end{aligned}\label{eq:simple iteration}
\end{equation}
for $\nu=0,1,2,\dots$, with an iteration parameter $\gamma\in(0,1)$.
Giving a prescribed tolerance and any initial guess of the density,
we can get a numerical approximation to the unknown density function
$\boldsymbol{\psi}$ once the simple iteration converges. However, 
it is difficult to calculate the boundary integral and volume integral
in Eqs. (\ref{eq:bie}) and (\ref{eq:simple iteration}) using the general
numerical integration methods since the analytic expressions of Green's 
functions, $i.e.$, the kernels of the integrals are unavailable. Hence, 
the KFBI method is a new strategy to evaluate the integrals approximately 
without using the analytical expressions of Green's functions which will
be described in the next section. The BIE (\ref{eq:bie}) can be further 
noted as $\mathcal{A}\boldsymbol{\psi}=\mathbf{g}$ 
where $\mathcal{A}\boldsymbol{\psi}\equiv\dfrac{1}{2}\boldsymbol{\psi}(\mathbf{p})-\int_{\partial\Omega}\mathbf{D}\mathbf{n}_{\mathbf{p}}\cdot\nabla_{\mathbf{p}}G(\mathbf{q};\mathbf{p})\boldsymbol{\psi}(\mathbf{q})ds_{\mathbf{q}}$
and $\mathbf{g}\equiv\mathbf{g}^{N}-\mathbf{n_{p}}\cdot\mathbf{D}\nabla_{\mathbf{p}}\left(\int_{\Omega}G(\mathbf{q};\mathbf{p})\mathbf{f}(\mathbf{q})d\mathbf{q}\right)$.
The simple iterations for solving the BIE $\mathcal{A}\boldsymbol{\psi}=\mathbf{g}$
could be further improved by a Krylov subspace method such as the GMRES iteration
\cite{saad1986gmres,saad2003iterative}.

\section{The Kernel-Free Boundary Integral Method \label{sec:KFBI}}

Once again, the KFBI method is used to approximately evaluate the volume  
and boundary integrals arising from the right hand side $\mathbf{g}$ of
the BIE $\mathcal{A}\boldsymbol{\psi}=\mathbf{g}$ and the matrix-vector 
multiplication $\mathcal{A}\boldsymbol{\psi}$ during the Krylov subspace 
iteration \cite{saad2003iterative}. The KFBI method first replaces the 
evaluation of the integrals by the solutions of the equivalent simple
interface problems on the larger regular domain $\mathcal{B}$ without calling 
for computing any approximation of Green's functions. After discretizing
the simple interface problem by the FDM on a Cartesian grid of the larger
regular domain $\mathcal{B}$, the obtained discrete system can be solved 
efficiently by a fast Fourier transform (FFT) based solver after applying 
an appropriate correction to the right hand side of the discrete system
at the irregular grid nodes. 
Noticing that the boundary and volume integrals in Eqs. (\ref{eq:bie}) 
and (\ref{eq:simple iteration}) are restricted at the discretization nodes 
of the interface $\partial\Omega$, we approximate the limiting values 
of the integrals by taking polynomial interpolation of the discrete 
numerical solutions. 

\subsection{Reinterpret the integrals\label{subsec:5.1}}

Assume that all of the variables and functions encountered are smooth 
enough such that the derivatives that appear are meaningful. A piecewise 
smooth dependent variable $\omega(\mathbf{p})$ which is defined on the 
larger regular domain $\mathcal{B}$ have possible discontinuities only 
on the domain boundary $\partial\Omega$. let $\omega^{+}(\mathbf{p})$ 
and $\omega^{-}(\mathbf{p})$ be the restrictions of $\omega(\mathbf{p})$ 
on the subdomains $\Omega$ and $\Omega^{C}$, respectively. For 
$\mathbf{p}\in\partial\Omega$, let 

\begin{equation}
\omega^{+}(\mathbf{p})\equiv\lim\limits_{
	\substack{\mathbf{q}\in\Omega\\
		\mathbf{q}\rightarrow \mathbf{p}}}\omega(\mathbf{q})
\quad\quad \text{and} \quad\quad 
\omega^{-}(\mathbf{p})\equiv\lim\limits_{
	\substack{\mathbf{q}\in\Omega^{C}\\
		\mathbf{q}\rightarrow \mathbf{p}}}\omega(\mathbf{q})
\end{equation}
be the limit values of $\omega(\mathbf{p})$ from either side of the domain
boundary. The jump of the variable $\omega(\mathbf{p})$ across the interface
is denoted by

\begin{equation}
[\omega(\mathbf{p})]\equiv \omega^{+}(\mathbf{p})-\omega^{-}(\mathbf{p})\label{eq:jump}
\end{equation}
with $\mathbf{p}\in\partial\Omega$. The jumps of a vector functions across the
interface can be defined in the same way. 

The volume integral $\mathbf{v}(\mathbf{p})=(\mathcal{G}\mathbf{f})(\mathbf{p})\equiv\int_{\Omega}G(\mathbf{q};\mathbf{p})\mathbf{f}(\mathbf{q})d\mathbf{q}$
satisfies the interface problem

\begin{equation}
\begin{array}{ccc}
\nabla\cdot(\mathbf{D}\nabla\mathbf{v})-\mathbf{K}\mathbf{v}=\widetilde{\mathbf{f}} &  & \textrm{in }\mathcal{B}\setminus\partial\Omega,\\{}
[\mathbf{v}]=0 &  & \textrm{on }\partial\Omega,\\
\mathbf{n}\cdot[\mathbf{D}\nabla\mathbf{v}]=0 &  & \textrm{on }\partial\Omega,\\
\mathbf{v}=0 &  & \textrm{on }\partial \mathcal{B}.
\end{array}\label{eq:volume}
\end{equation}
The equivalent interface problem of the volume integral is a partial
differential equation with Dirichlet boundary condition which can
be solved easily. Here, the source function in Eq.(\ref{eq:bvp}) $\mathbf{f}=(-\kappa V_{\text{m}}^{n},\kappa V_{\text{m}}^{n}-\beta I_{\text{stim}})^{T}$
is a vector function without analytical expression consisting of 
the discrete solution of the transmembrane voltage at the last 
discrete time. If we set the $\mathbf{\widetilde{f}}$ is zero in the
subdomain $\Omega^{C}$and equal to $\mathbf{f}$ in the subdomain
$\Omega$, it is not easy to compute the jump of $\mathbf{\widetilde{f}}$
at the discretization nodes of the domain boundary during solving 
the equivalent interface problems. According to the continuous 
properties of the volume integral and the boundary integral 
\cite{xie2019fourth,ying2007kernel}, we continuously extend the 
vector function $\mathbf{f}(\mathbf{p})$ ($\mathbf{p}\in\Omega$)
by evaluating the transmembrane voltage $V_{\text{m}}$ in the 
larger regular domain $\mathcal{B}$ at the very beginning and 
approximately estimate the value of 
$\mathbf{n_{p}}\cdot\mathbf{D}\nabla_{\mathbf{p}}\left(\int_{\Omega}G(\mathbf{q};\mathbf{p})\mathbf{f}(\mathbf{q})d\mathbf{q}\right)$ 
by taking polynomial interpolation of the discrete numerical solution of the 
equivalent interface problem, $i.e.$, the value of $\mathbf{g}$ in the 
BIE $\mathcal{A}\boldsymbol{\psi}=\mathbf{g}$.

The boundary integral $\mathbf{v}(\mathbf{p})=-(\mathcal{K}\boldsymbol{\psi})(\mathbf{p})\equiv-\int_{\partial\Omega}G(\mathbf{q};\mathbf{p})\boldsymbol{\psi}(\mathbf{q})ds_{\mathbf{q}}$
satisfies the interface problem

\begin{equation}
\begin{array}{ccc}
\nabla\cdot(\mathbf{D}\nabla\mathbf{v})-\mathbf{K}\mathbf{v}=0 &  & \textrm{in }\mathcal{B}\setminus\partial\Omega,\\{}
[\mathbf{v}]=0 &  & \textrm{on }\partial\Omega,\\
\mathbf{n}\cdot[\mathbf{D}\nabla\mathbf{v}]=\boldsymbol{\psi} &  & \textrm{on }\partial\Omega,\\
\mathbf{v}=0 &  & \textrm{on }\partial \mathcal{B}.
\end{array}\label{eq:single}
\end{equation}
The boundary integral is continuous across the interface $\partial\Omega$
while the normal flux of the boundary integral 
$\mathbf{n}_{\mathbf{p}}\cdot[\mathbf{D}\nabla_{\mathbf{p}}\mathbf{v}]$
has a jump \cite{ying2007kernel} with strength equal to $\boldsymbol{\psi}(\mathbf{p})$
, $i.e.$,

\begin{equation}
\begin{array}{ccc}
\mathbf{v}^{+}(\mathbf{p})=\mathbf{v}^{-}(\mathbf{p}) &  & \textrm{on}\;\partial\Omega\end{array}
\end{equation}
and

\begin{equation}
\begin{array}{ccc}
\mathbf{n}\cdot(\mathbf{D}\nabla\mathbf{v}^{+})(\mathbf{p})-\mathbf{n}\cdot(\mathbf{D}\nabla\mathbf{v})(\mathbf{p})=\dfrac{1}{2}\boldsymbol{\psi}(\mathbf{p}) &  & \textrm{on}\;\partial\Omega\end{array}\label{eq:single-use}
\end{equation}

\begin{equation}
\begin{array}{ccc}
\mathbf{n}\cdot(\mathbf{D}\nabla\mathbf{v}^{-})(\mathbf{p})-\mathbf{n}\cdot(\mathbf{D}\nabla\mathbf{v})(\mathbf{p})=-\dfrac{1}{2}\boldsymbol{\psi}(\mathbf{p}) &  & \textrm{on}\;\partial\Omega\end{array}
\end{equation}
For the interior boundary value problems, the discontinuity property
of the boundary integral presented in Eq.(\ref{eq:single-use}) is used to evaluate 
the left hand side of the BIE (\ref{eq:bie}), $i.e.$, the matrix-vector
multiplication $\mathcal{A}\boldsymbol{\psi}$. Moreover, the solutions
of the above interface problems are both continuous across the interface
which makes the extension of the $\mathbf{f}$ to $\widetilde{\mathbf{f}}$ 
be reasonable. We point out that the above two interface problems defined in
Eqs. (\ref{eq:volume}) and (\ref{eq:single}) can be presented and solved 
in a unified framework.

\subsection{Discretization of the PDE on a Cartesian grid }

We discretize the partial differential equation 

\begin{equation}
\begin{array}{cc}
\mathcal{L}\mathbf{u}(\mathbf{p})\equiv\nabla\cdot(\mathbf{D}\nabla\mathbf{u})-\mathbf{K}\mathbf{u}=\mathbf{f}(\mathbf{p}), & \mathbf{p}\in\mathcal{B}\setminus\partial\Omega\end{array}\label{eq:PDE}
\end{equation}
on a Cartesian grid of the larger regular domain $\mathcal{B}$ with the FDM. 
In this work, we only present the finite difference discretization scheme in 
two-dimensional space since the finite difference discretization scheme 
in three-dimensional space is similar to that in two-dimensional space.
For simplicity, we assume that the larger regular domain $\mathcal{B}=(a,b)\times(c,d)$ 
is a square which is partitioned into a uniform $I\times I$ Cartesian grid 
with the mesh parameter $h=(b-a)/I$. Let $x_{k}=a+kh$ and $y_{l}=c+lh$ for 
$k=0,1,\cdots,I$ and $l=0,1,\cdots,I$. Denote the $(k,l)$-th node of the 
Cartesian grid by $\mathbf{p}_{k,l}=(x_{k},y_{l})^{T}$. Omitting the 
superscripts representing the discrete times, we denote the finite difference 
approximation of the solution $\mathbf{u}$ at the $(k,l)$-th grid node 
$\mathbf{p}_{k,l}$ by $\mathbf{u}_{k,l}=(\Phi_{k,l}^{\text{i}},\Phi_{k,l}^{\text{e}})^{T}$.
For the bidomain simulations, we model the square as consisting of fibers 
that lie parallel to the $x$-axis, $i.e.$, the conductivity tensors 
$\mathbf{D}_{\text{i}}$ and $\mathbf{D}_e$ are diagonal matrices. Denote 
the intra- and extracellular conductivity tensors by
$\mathbf{D}_{\text{i}}=\textrm{diag}(\sigma_{\text{i}},\eta_{\text{i}})$ and
$\mathbf{D}_{\text{e}}=\textrm{diag}(\sigma_{\text{e}},\eta_{\text{e}})$ and
denote the source function by $\mathbf{f}=(f^{\text{i}},f^{\text{e}})^{T}$. 
With the second-order centered finite difference method, the PDE (\ref{eq:PDE}) 
is discretized into the finite difference equations

\begin{equation}
\begin{aligned}\dfrac{s_{k,l}-2\sigma_{\text{i}}\Phi_{k,l}^{\text{i}}-2\eta_{\text{i}}\Phi_{k,l}^{\text{i}}}{h^{2}}-\kappa\Phi_{k,l}^{\text{i}}+\kappa\Phi_{k,l}^{\text{e}} & =f_{k,l}^{\text{i}},\\
\dfrac{t_{k,l}-2\sigma_{\textrm{e}}\Phi_{k,l}^{\text{e}}-2\eta_{\text{e}}\Phi_{k,l}^{\text{e}}}{h^{2}}+\kappa\Phi_{k,l}^{\text{i}}-\kappa\Phi_{k,l}^{\text{e}} & =f_{k,l}^{\text{e}},
\end{aligned}
\label{eq:FDE}
\end{equation}
where 

\begin{equation}
s_{k,l}=\sigma_{\text{i}}(\Phi_{k+1,l}^{\text{i}}+\Phi_{k-1,l}^{\text{i}})+\eta_{\text{i}}(\Phi_{k,l+1}^{\text{i}}+\Phi_{k,l-1}^{\text{i}}),
\end{equation}

\begin{equation}
t_{k,l}=\sigma_{\text{e}}(\Phi_{k+1,l}^{\text{e}}+\Phi_{k-1,l}^{\text{e}})+\eta_{\text{e}}(\Phi_{k,l+1}^{\text{e}}+\Phi_{k,l-1}^{\text{e}}),
\end{equation}
for $k=0,1,\dots,I$ and $l=0,1,\dots,I$. In two- and three-dimensional spaces, 
the finite difference equation at each Cartesian grid node $\mathbf{p}_{k,l}$ 
involves a five-point stencil and a seven-point stencil, respectively. When
the coefficients $\sigma_{\text{i}}$, $\sigma_{\text{e}}$, $\eta_{\text{i}}$, 
$\eta_{\text{e}}$, $\kappa$, and the source function $\mathbf{f}$ are 
sufficiently smooth, the solution to the finite difference equation (\ref{eq:FDE})
has a second-order accuracy, irrespective of the interface $\partial\Omega$.

\subsection{Correction of the discrete system }

Owing to the existence of the interface $\partial\Omega$, the Cartesian grid 
nodes $\mathbf{p}_{k,l}$ are divided into regular nodes and irregular
nodes. A grid node is called irregular if its finite difference stencil 
has intersections with the domain boundary $\partial\Omega$. Otherwise, the grid 
nodes are called regular. The local truncation
errors of the finite difference equations (\ref{eq:FDE}) at irregular
grid nodes are much larger than those at regular grid nodes. Specifically,
the estimate for the local truncation errors at irregular grid nodes
is $\mathcal{O}(h^{-2})$ with mesh parameter $h$ \cite{ying2014kernel}.
Because of the large local truncation errors at irregular grid nodes, the 
solution to the discrete system obtained from the equivalent interface 
problem is inaccurate. Hence, we need to correct the finite difference 
equations at irregular grid nodes to get an accurate solution.

Denote the smooth solution in subdomains $\Omega$ and $\Omega^{C}$
by $\mathbf{u}^{+}$ and $\mathbf{u}^{-}$. In the case of two-dimensional 
space, we assume that the five-point stencil of an irregular grid node
$\mathbf{p}_{k,l}$ intersects the domain boundary $\partial\Omega$ only
once and the intersection point located at the horizontal grid line is 
$(z_{k},y_{l})$ with $x_{k}\leq z_{k}<x_{k+1}$. In the next, we estimate
the local truncation error at the irregular grid node 
$\mathbf{p}_{k,l}$ 

\begin{equation}
E_{h,x}^{\text{i}}(x_{k},y_{l})\equiv\dfrac{\sigma_{\text{i}}\left((\Phi^{\text{i}})^{-}(x_{k+1},y_{l})+(\Phi^{\text{i}})^{+}(x_{k-1},y_{l})-2(\Phi^{\text{i}})^{+}(x_{k},y_{l})\right)}{h^{2}}-\sigma_{\text{i}}\dfrac{\partial^{2}}{\partial x^{2}}(\Phi^{\text{i}})^{+}(x_{k},y_{l}),
\end{equation}

\begin{equation}
E_{h,x}^{^{\text{e}}}(x_{k},y_{l})\equiv\dfrac{\sigma_{\text{e}}\left((\Phi^{\text{e}})^{-}(x_{k+1},y_{l})+(\Phi^{\text{e}})^{+}(x_{k-1},y_{l})-2(\Phi^{\text{e}})^{+}(x_{k},y_{l})\right)}{h^{2}}-\sigma_{\text{e}}\dfrac{\partial^{2}}{\partial x^{2}}(\Phi^{\text{e}})^{+}(x_{k},y_{l}).
\end{equation}
We first make Taylor expansions for $(\Phi^{\text{i}})^{+}(x_{k+1},y_{l})$, 
$(\Phi^{\text{e}})^{+}(x_{k+1},y_{l})$, $(\Phi^{\text{i}})^{-}(x_{k+1},y_{l})$, 
and $(\Phi^{\text{e}})^{-}(x_{k+1},y_{l})$ around the intersection point 
$(z_{k},y_{l})$,

\begin{equation*}
\begin{aligned}
(\Phi^{\textrm{i}})^{\pm}(x_{k+1},y_{l}) = &(\Phi^{\textrm{i}})^{\pm}(z_{k},y_{l})+\partial_{x}(\Phi^{\textrm{i}})^{\pm}(z_{k},y_{l})(x_{k+1}-z_{k})+\dfrac{1}{2}\partial_{xx}(\Phi^{\textrm{i}})^{\pm}(z_{k},y_{l})(x_{k+1}-z_{k})^{2}\\
& + \dfrac{1}{6}\partial_{xxx}(\Phi^{\textrm{i}})^{\pm}(z_{k},y_{l})(x_{k+1}-z_{k})^{3}+O(h^{4}),
\end{aligned}
\end{equation*}
\noindent and 

\begin{equation*}
\begin{aligned}
(\Phi^{\textrm{e}})^{\pm}(x_{k+1},y_{l}) = &(\Phi^{\textrm{e}})^{\pm}(z_{k},y_{l})+\partial_{x}(\Phi^{\textrm{e}})^{\pm}(z_{k},y_{l})(x_{k+1}-z_{k})+\dfrac{1}{2}\partial_{xx}(\Phi^{\textrm{e}})^{\pm}(z_{k},y_{l})(x_{k+1}-z_{k})^{2}\\
& + \dfrac{1}{6}\partial_{xxx}(\Phi^{\textrm{e}})^{\pm}(z_{k},y_{l})(x_{k+1}-z_{k})^{3}+O(h^{4}).
\end{aligned}
\end{equation*}
We then estimate the leading order term of the local truncation error 
$E_{h,x}^{^{\textrm{i}}}(x_{k},y_{l})$ and $E_{h,x}^{^{\textrm{e}}}(x_{k},y_{l})$ as 

\begin{equation*}
\begin{aligned}
\dfrac{(\Phi^{\textrm{i}})^{+}(x_{k+1},y_{l})-(\Phi^{\textrm{i}})^{-}(x_{k+1},y_{l})}{h^{2}}=&\dfrac{1}{h^{2}}\left\{ [\Phi^{\textrm{i}}]+[\partial_{x}\Phi^{\textrm{i}}](x_{k+1}-z_{k})+\dfrac{1}{2}[\partial_{xx}\Phi^{\textrm{i}}](x_{k+1}-z_{k})^{2}\right. \\
&\left. +\dfrac{1}{6}[\partial_{xxx}\Phi^{\textrm{i}}](x_{k+1}-z_{k})^{3}\right\} +O(h^{2})
\end{aligned}
\end{equation*}
\noindent and 

\begin{equation*}
\begin{aligned}
\dfrac{(\Phi^{\textrm{e}})^{+}(x_{k+1},y_{l})-(\Phi^{\textrm{e}})^{-}(x_{k+1},y_{l})}{h^{2}}=&\dfrac{1}{h^{2}}\left\{ [\Phi^{\textrm{e}}]+[\partial_{x}\Phi^{\textrm{e}}](x_{k+1}-z_{k})+\dfrac{1}{2}[\partial_{xx}\Phi^{\textrm{e}}](x_{k+1}-z_{k})^{2} \right.\\
&\left. +\dfrac{1}{6}[\partial_{xxx}\Phi^{\textrm{e}}](x_{k+1}-z_{k})^{3}\right\} +O(h^{2}).
\end{aligned}
\end{equation*}
Here, $[\Phi^{\text{i}}]=(\Phi^{\text{i}})^{+}(z_{k},y_{l})-(\Phi^{\text{i}})^{-}(z_{k},y_{l})$,
$[\Phi^{\text{e}}]=(\Phi^{\text{e}})^{+}(z_{k},y_{l})-(\Phi^{\text{e}})^{-}(z_{k},y_{l})$,
$[\partial_{x}\Phi^{\text{i}}]=\partial_{x}(\Phi^{\text{i}})^{+}(z_{k},y_{l})-\partial_{x}(\Phi^{\text{i}})^{-}(z_{k},y_{l})$, $[\partial_{x}\Phi^{\text{e}}]=\partial_{x}(\Phi^{\text{e}})^{+}(z_{k},y_{l})-\partial_{x}(\Phi^{\text{e}})^{-}(z_{k},y_{l})$, 
$[\partial_{xx}\Phi^{\text{i}}]\equiv\partial_{xx}(\Phi^{\text{i}})^{+}(z_{k},y_{l})-\partial_{xx}(\Phi^{\text{i}})^{-}(z_{k},y_{l})$,
and $[\partial_{xx}\Phi^{\text{e}}]\equiv\partial_{xx}(\Phi^{\text{e}})^{+}(z_{k},y_{l})-\partial_{xx}(\Phi^{\text{e}})^{-}(z_{k},y_{l})$ are the jumps of the function and its partial derivatives across the
interface. 

\begin{figure}[!t]
	\centering
	\includegraphics[scale=.5]{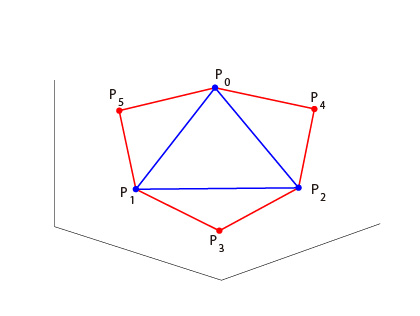}
	\caption{Six boundary points of a local coordinate system.\label{fig:sixNeighbour}}
\end{figure}

The correction of the discrete system at the irregular grid node $\mathbf{p}_{k,l}$
is given by the following quantity
\begin{equation}
(C_{h,x}^{\text{i}})^{+}=\dfrac{\sigma_{\text{i}}}{h^{2}}\left([\Phi^{\text{i}}]+[\partial_{x}\Phi^{\text{i}}](x_{k+1}-z_{k})+\dfrac{1}{2}[\partial_{xx}\Phi^{\text{i}}](x_{k+1}-z_{k})^{2}\right)
\end{equation}

\begin{equation}
(C_{h,x}^{\text{e}})^{+}=\dfrac{\sigma_{\text{e}}}{h^{2}}\left([\Phi^{\text{e}}]+[\partial_{x}\Phi^{\text{e}}](x_{k+1}-z_{k})+\dfrac{1}{2}[\partial_{xx}\Phi^{\text{e}}](x_{k+1}-z_{k})^{2}\right)
\end{equation}
The jumps of the function and its partial derivatives across the interface
are derived from the interface conditions of the equivalent interface problem 
\cite{ying2007kernel,ying2014kernel}. In three-dimensional space, we calculate
the corresponding jumps by constructing local coordinate systems on the 
discrete triangle grid which represents the domain boundary $\partial \Omega$.
As shown in Fig. \ref{fig:sixNeighbour}, the local coordinate system of 
a triangle element consists of the six boundary points around it. 

When the interface intersects the finite difference stencil multiple times, 
we should calculate all the corresponding terms $(C_{h,x}^{\text{i}})^{\pm}$, 
$(C_{h,x}^{\text{e}})^{\pm}$, $(C_{h,y}^{\text{i}})^{\pm}$, $(C_{h,y}^{\text{e}})^{\pm}$,
$(C_{h,z}^{\text{i}})^{\pm},$ and $(C_{h,z}^{\text{e}})^{\pm}$ to modify the 
right hand side of the discrete linear system (\ref{eq:PDE}).
In general, for an irregular grid node $(x_{k},y_{l})$ in two-dimensional space, 
the finite difference equation turns to 

\begin{equation}
\begin{aligned}\dfrac{s_{k,l}-2\sigma_{\text{i}}\Phi_{k,l}^{\text{i}}-2\eta_{\text{i}}\Phi_{k,l}^{\text{i}}}{h^{2}}-\kappa\Phi_{k,l}^{\text{i}}+\kappa\Phi_{k,l}^{\text{e}} & =f_{k,l}^{\text{i}}+C_{h}^{\text{i}},\\
\dfrac{t_{k,l}}{h^{2}-2\sigma_{\text{e}}\Phi_{k,l}^{\text{e}}-2\eta_{\text{e}}\Phi_{k,l}^{\text{e}}}+\kappa\Phi_{k,l}^{\text{i}}-\kappa\Phi_{k,l}^{\text{e}} & =f_{k,l}^{\text{e}}+C_{h}^{\text{e}}.
\end{aligned}\label{eq:correct}
\end{equation}
where $C_{\text{h}}^{\text{i}}$ and $C_{\text{h}}^{\text{e}}$ are the sum of 
all corrections at the irregular grid node $(x_{k},y_{l})$ of the variable 
$\Phi^{\text{i}}$ and $\Phi^{\text{e}}$, respectively. 
After appropriately modifying the right hand side of the discrete linear
system, the local truncation error of the obtained discrete linear
system is of first-order at irregular grid nodes and its solution
has essentially second-order accuracy \cite{beale2007accuracy,ying2007kernel}. 
Since the modification at irregular grid nodes does not change the
symmetric and positive definite coefficient matrix, the discrete 
linear system of the equivalent interface problem can be solved 
efficiently by an FFT based Poisson solver 
\cite{dorr1970direct,swarztrauber1977methods}. 

\subsection{Interpolation of the volume and boundary integrals\label{subsec:5.4}}

\begin{figure}[!t]
	\centering
	\includegraphics[scale=0.5]{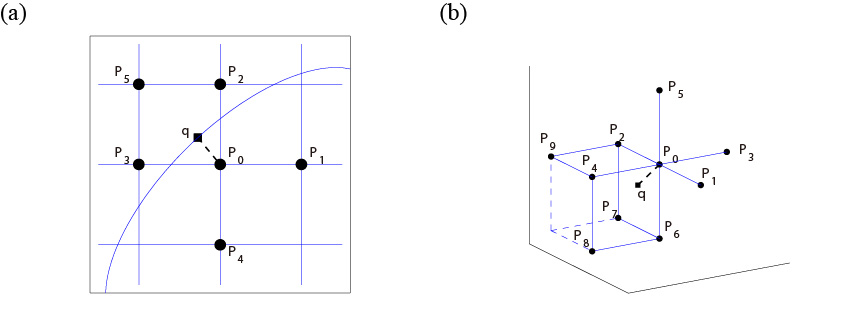}
	\caption{Illustration of a six-point and a ten-point interpolation stencil of
		the discretization node $\mathbf{q}$ of the boundary curve.\label{fig:six_stencil}}	
\end{figure}

Denote $\mathbf{u}_{\nu,h}(\mathbf{p})=(\Phi_{\nu,h}^{\textrm{i}},\Phi_{\nu,h}^{\textrm{e}})^T$ 
as the solution of the discrete 
linear system obtained after appropriate correction at the irregular 
grid nodes and $\nu$ is the iteration number of the iteration method 
for solving the BIE (\ref{eq:bie}). We make Taylor expansions for the 
approximate solution $\mathbf{u}_{\nu,h}$ is to extract the limit values 
of $\mathbf{n}\cdot(\mathbf{D}\nabla\mathbf{u}_{\nu,h}^{+})$
, $i.e.$, the values of
$\partial_{x}\Phi_{\text{i}}^{+}$, $\partial_{x}\Phi_{\text{e}}^{+}$, 
$\partial_{y}\Phi_{\text{i}}^{+}$, $\partial_{y}\Phi_{\text{e}}^{+}$, 
$\partial_{z}\Phi_{\text{i}}^{+}$, and $\partial_{z}\Phi_{\text{e}}^{+}$ 
at the discretization nodes of the boundary curve $\partial \Omega$.
As shown in Fig. \ref{fig:six_stencil}, a six-point interpolation stencil
in two-dimensional space or a ten-point interpolation stencil in 
three-dimensional space is required to evaluate the limit values
at a discretization node $\mathbf{q}$ of the boundary curve by 
quadratic interpolation. 

For a Cartesian grid node $\mathbf{p}\in\Omega\subset\mathbb{R}^{2}$, the Taylor
expansion for the approximate solution $\mathbf{u}_{\nu,h}(\mathbf{p})=(\Phi_{\nu,h}^{\text{i}}(\mathbf{p}),\Phi_{\nu,h}^{\text{e}}(\mathbf{p}))^{T}$
around the discretization node $\mathbf{q}$ of the boundary curve is as followed

\begin{equation}
\begin{aligned}\Phi_{\nu,h}^{\text{i}}(\mathbf{p}) & =(\Phi_{\nu,h}^{\text{i}})^{+}(\mathbf{q})+\dfrac{\partial(\Phi_{\nu,h}^{\text{i}})^{+}(\mathbf{q})}{\partial x}\xi+\dfrac{\partial(\Phi_{\nu,h}^{\text{i}})^{+}(\mathbf{q})}{\partial y}\eta\\
& +\dfrac{1}{2}\dfrac{\partial^{2}(\Phi_{\nu,h}^{\text{i}})^{+}(\mathbf{q})}{\partial x^{2}}\xi^{2}+\dfrac{\partial^{2}(\Phi_{\nu,h}^{\text{i}})^{+}(\mathbf{q})}{\partial x\partial y}\xi\eta+\dfrac{1}{2}\dfrac{\partial^{2}(\Phi_{\nu,h}^{\text{i}})^{+}(\mathbf{q})}{\partial y^{2}}\eta^{2}+O(|\mathbf{p}-\mathbf{q}|^{3}),
\end{aligned}
\end{equation}

\begin{equation}
\begin{aligned}\Phi_{\nu,h}^{\text{e}}(\mathbf{p}) & =(\Phi_{\nu,h}^{\text{e}})^{+}(\mathbf{q})+\dfrac{\partial(\Phi_{\nu,h}^{\text{e}})^{+}(\mathbf{q})}{\partial x}\xi+\dfrac{\partial(\Phi_{\nu,h}^{\text{e}})^{+}(\mathbf{q})}{\partial y}\eta\\
& +\dfrac{1}{2}\dfrac{\partial^{2}(\Phi_{\nu,h}^{\text{e}})^{+}(\mathbf{q})}{\partial x^{2}}\xi^{2}+\dfrac{\partial^{2}(\Phi_{\nu,h}^{\text{e}})^{+}(\mathbf{q})}{\partial x\partial y}\xi\eta+\dfrac{1}{2}\dfrac{\partial^{2}(\Phi_{\nu,h}^{\text{e}})^{+}(\mathbf{q})}{\partial y^{2}}\eta^{2}+O(|\mathbf{p}-\mathbf{q}|^{3}).
\end{aligned}
\end{equation}
For a Cartesian grid node $\mathbf{p}\in\Omega^{C}\subset\mathbb{R}^{2}$, the Taylor expansion of
the approximate solution $\mathbf{u}_{\nu,h}(\mathbf{p})=(\Phi_{\nu,h}^{\text{i}}(\mathbf{p}),\Phi_{\nu,h}^{\text{e}}(\mathbf{p}))^{T}$
around the discretization node $\mathbf{q}$ of the boundary curve is as followed

\begin{equation}
\begin{aligned}\Phi_{\nu,h}^{\text{i}}(\mathbf{p}) & =(\Phi_{\nu,h}^{\text{i}})^{-}(\mathbf{q})+\dfrac{\partial(\Phi_{\nu,h}^{\text{i}})^{-}(\mathbf{q})}{\partial x}\xi+\dfrac{\partial(\Phi_{\nu,h}^{\text{i}})^{-}(\mathbf{q})}{\partial y}\eta\\
& +\dfrac{1}{2}\dfrac{\partial^{2}(\Phi_{\nu,h}^{\text{i}})^{-}(\mathbf{q})}{\partial x^{2}}\xi^{2}+\dfrac{\partial^{2}(\Phi_{\nu,h}^{\text{i}})^{-}(\mathbf{q})}{\partial x\partial y}\xi\eta+\dfrac{1}{2}\dfrac{\partial^{2}(\Phi_{\nu,h}^{\text{i}})^{-}(\mathbf{q})}{\partial y^{2}}\eta^{2}+O(|\mathbf{p}-\mathbf{q}|^{3}),
\end{aligned}
\end{equation}

\begin{equation}
\begin{aligned}\Phi_{\nu,h}^{\text{e}}(\mathbf{p}) & =(\Phi_{\nu,h}^{\text{e}})^{-}(\mathbf{q})+\dfrac{\partial(\Phi_{\nu,h}^{\text{e}})^{-}(\mathbf{q})}{\partial x}\xi+\dfrac{\partial(\Phi_{\nu,h}^{\text{e}})^{-}(\mathbf{q})}{\partial y}\eta\\
& +\dfrac{1}{2}\dfrac{\partial^{2}(\Phi_{\nu,h}^{\text{e}})^{-}(\mathbf{q})}{\partial x^{2}}\xi^{2}+\dfrac{\partial^{2}(\Phi_{\nu,h}^{\text{e}})^{-}(\mathbf{q})}{\partial x\partial y}\xi\eta+\dfrac{1}{2}\dfrac{\partial^{2}(\Phi_{\nu,h}^{\text{e}})^{-}(\mathbf{q})}{\partial y^{2}}\eta^{2}+O(|\mathbf{p}-\mathbf{q}|^{3}).
\end{aligned}
\end{equation}
Here, $(\xi,\eta)^{T}\equiv\mathbf{p}-\mathbf{q}$. Omit some subscripts
by

\begin{equation}
\begin{array}{ccc}
V^{\pm}\equiv(\Phi_{\nu,h}^{\text{i}})^{\pm}(\mathbf{q}), & V_{x}^{\pm}\equiv\dfrac{\partial(\Phi_{\nu,h}^{\text{i}})^{\pm}(\mathbf{q})}{\partial x}, & V_{y}^{\pm}\equiv\dfrac{\partial(\Phi_{\nu,h}^{\text{i}})^{\pm}(\mathbf{q})}{\partial y},\\
V_{xx}^{\pm}\equiv\dfrac{\partial^{2}(\Phi_{\nu,h}^{\text{i}})^{\pm}(\mathbf{q})}{\partial x^{2}}, & V_{xy}^{\pm}\equiv\dfrac{\partial^{2}(\Phi_{\nu,h}^{\text{i}})^{\pm}(\mathbf{q})}{\partial x\partial y}, & V_{yy}^{\pm}\equiv\dfrac{\partial^{2}(\Phi_{\nu,h}^{\text{i}})^{\pm}(\mathbf{q})}{\partial y^{2}},
\end{array}
\end{equation}
and 

\begin{equation}
\begin{array}{ccc}
W^{\pm}\equiv(\Phi_{\nu,h}^{\text{e}})^{\pm}(\mathbf{q}), & W_{x}^{\pm}\equiv\dfrac{\partial(\Phi_{\nu,h}^{\text{e}})^{\pm}(\mathbf{q})}{\partial x}, & W_{y}^{\pm}\equiv\dfrac{\partial(\Phi_{\nu,h}^{\text{e}})^{\pm}(\mathbf{q})}{\partial y},\\
W_{xx}^{\pm}\equiv\dfrac{\partial^{2}(\Phi_{\nu,h}^{\text{e}})^{\pm}(\mathbf{q})}{\partial x^{2}}, & W_{xy}^{\pm}\equiv\dfrac{\partial^{2}(\Phi_{\nu,h}^{\text{e}})^{\pm}(\mathbf{q})}{\partial x\partial y}, & W_{yy}^{\pm}\equiv\dfrac{\partial^{2}(\Phi_{\nu,h}^{\text{e}})^{\pm}(\mathbf{q})}{\partial y^{2}}.
\end{array}
\end{equation}
Denote the six nearly grid nodes by $\mathbf{p}_{j}$ $(j=0,1,\cdots,5)$. First, 
we evaluate the truncated Taylor series of $\Phi_{\nu,h}^{\text{i}}$ which yield

\begin{equation}
\begin{aligned}V^{+}+V_{x}^{+}\xi_{j}+V_{y}^{+}\eta_{j}+\dfrac{1}{2}V_{xx}^{+}\xi_{j}^{2}+V_{xy}^{+}\xi_{j}\eta_{j}+\dfrac{1}{2}V_{yy}^{+}\eta_{j}^{2}=V_{j} & \;\textrm{ if}\;\mathbf{P}_{j}\in\Omega\end{aligned}
\label{eq:taylor-i-po}
\end{equation}
and 
\begin{equation}
\begin{aligned}V^{-}+V_{x}^{-}\xi_{j}+V_{y}^{-}\eta_{j}+\dfrac{1}{2}V_{xx}^{-}\xi_{j}^{2}+V_{xy}^{-}\xi_{j}\eta_{j}+\dfrac{1}{2}V_{yy}^{-}\eta_{j}^{2}=V_{j} & \;\textrm{ if}\;\mathbf{P}_{j}\in\Omega^{C},\end{aligned}
\label{eq:taylor-i}
\end{equation}
where $V_{j}\equiv(\Phi_{\nu,h}^{\text{i}})(\mathbf{p}_{j})$ and
$(\xi_{j},\eta_{j})^{T}\equiv(\mathbf{p}_{j}-\mathbf{q})$
for $j=0,1,\cdots,5$. Similarly, we evaluate
the truncated Taylor series of $\Phi_{\nu,h}^{\text{e}}$ which yield 

\begin{equation}
\begin{aligned}W^{+}+W_{x}^{+}\xi_{j}+W_{y}^{+}\eta_{j}+\dfrac{1}{2}W_{xx}^{+}\xi_{j}^{2}+W_{xy}^{+}\xi_{j}\eta_{j}+\dfrac{1}{2}W_{yy}^{+}\eta_{j}^{2}=W_{j} & \;\textrm{ if}\;\mathbf{P}_{j}\in\Omega\end{aligned}
\label{eq:taylor-e-po}
\end{equation}
and 
\begin{equation}
\begin{aligned}W^{-}+W_{x}^{-}\xi_{j}+W_{y}^{-}\eta_{j}+\dfrac{1}{2}W_{xx}^{-}\xi_{j}^{2}+W_{xy}^{-}\xi_{j}\eta_{j}+\dfrac{1}{2}W_{yy}^{-}\eta_{j}^{2}=W_{j} & \;\textrm{ if}\;\mathbf{P}_{j}\in\Omega^{C},\end{aligned}
\label{eq:taylor-e}
\end{equation}
where $W_{j}\equiv(\Phi_{\nu,h}^{\text{e}})(\mathbf{p}_{j})$ for $j=0,1,\cdots,5$.
Let 

\begin{equation}
J_{j}^{\text{i}}\equiv[V]+[V_{x}]\xi_{j}+[V_{y}]\eta_{j}+\dfrac{1}{2}[V_{xx}]\xi_{j}^{2}+[V_{xy}]\xi_{j}\eta_{j}+\dfrac{1}{2}[V_{yy}]\eta_{j}^{2}\label{eq:ji}
\end{equation}
and 

\begin{equation}
J_{j}^{\text{e}}=[W]+[W_{x}]\xi_{j}+[W_{y}]\eta_{j}+\dfrac{1}{2}[W_{xx}]\xi_{j}^{2}+[W_{xy}]\xi_{j}\eta_{j}+\dfrac{1}{2}[W_{yy}]\eta_{j}^{2},\label{eq:je}
\end{equation}
for $j=0,1,\cdots,5$. If the grid node $p_{j}\in\Omega^{C}$, we rewrite
Eq. (\ref{eq:taylor-i}) and Eq. (\ref{eq:taylor-e}) as

\begin{equation}
\begin{aligned}V^{+}+V_{x}^{+}\xi_{j}+V_{y}^{+}\eta_{j}+\dfrac{1}{2}V_{xx}^{+}\xi_{j}^{2}+V_{xy}^{+}\xi_{j}\eta_{j}+\dfrac{1}{2}V_{yy}^{+}\eta_{j}^{2}=V_{j}+J_{j}^{\text{i}} & \;\textrm{ if}\;\mathbf{P}_{j}\in\Omega^{C}\end{aligned}
\label{eq:taylor-i-co}
\end{equation}
and 

\begin{equation}
\begin{aligned}W^{+}+W_{x}^{+}\xi_{j}+W_{y}^{+}\eta_{j}+\dfrac{1}{2}W_{xx}^{+}\xi_{j}^{2}+W_{xy}^{+}\xi_{j}\eta_{j}+\dfrac{1}{2}W_{yy}^{+}\eta_{j}^{2}=W_{j}+J_{j}^{\text{e}} & \;\textrm{ if}\;\mathbf{P}_{j}\in\Omega^{C}.\end{aligned}
\label{eq:taylor-e-co}
\end{equation}

We need to calculate the jumps encountered in Eqs. (\ref{eq:ji}) and (\ref{eq:je})
at the discretization nodes of the boundary curve $\partial\Omega$. 
Using the values of $(\xi_{j},\eta_{j})^{T}$ and $(V_{j},W_{j})^{T}$ for $j=0,1,\cdots,5$,
we solve the linear system in Eqs. (\ref{eq:taylor-i-po}) and (\ref{eq:taylor-i-co})
for $V^{+}$, $V_{x}^{+}$, $V_{y}^{+}$, $V_{xx}^{+}$, $V_{xy}^{+}$,
$V_{yy}^{+}$ and the linear system in Eqs. (\ref{eq:taylor-e-po})
and (\ref{eq:taylor-e-co}) for $W^{+}$, $W_{x}^{+}$, $W_{y}^{+}$,
$W_{xx}^{+}$, $W_{xy}^{+}$, $W_{yy}^{+}$, respectively. In the implementation,
we usually rewrite the coefficient matrix as independent of the mesh
parameter $h$ by changing $\xi_{j}$ to $\xi_{j}/h$ and changing
$\eta_{j}$ to $\eta_{j}/h$ for $j=0,1,\cdots,5$. In the case of the 
two-dimensional space, if we choose the six grid nodes $\mathbf{P}_{j}$ $(j=0,1,\dots,5)$ 
appropriately, $e.g.$, the nearest six grid nodes in a structured grid, 
the limit values of the solution $\mathbf{u}_{\nu,h}$ and its derivatives 
are uniquely determined. As shown in Fig. \ref{fig:six_stencil}, 
we choose the nearest ten grid nodes in the structured grid as 
a ten-point interpolation stencil in the three-dimensional space. 

\section{Algorithm Summary}

In this section, we summary the algorithm with emphasis on numerically 
solving the diffusion part of the bidomain equations by the KFBI method. 
After choosing the appropriate operator splitting technique and temporal 
discretization scheme, we present the procedure of solving the BVP 
defined in Eq.(\ref{eq:pde1}) or Eq.(\ref{eq:bvp}) with homogeneous 
Neumann boundary condition by the simple Richardson iteration \ref{alg:algo}
from $t^n$ to $t^{n+1}$.

\begin{algorithm}[H] \label{alg:algo}	
	\caption{the KFBI method}
	Denote . Give an initial guess as $\boldsymbol{\psi}_{\nu}=0$ $(\nu=0)$
	where $\nu$ is the iteration number of the simple Richardson iteration.\\
	Denote $\mathbf{r}=\mathbf{g}^{N}-\mathbf{n_{p}}\cdot\mathbf{D}\nabla_{\mathbf{p}}\left(\int_{\Omega}G(\mathbf{q};\mathbf{p})\mathbf{f}(\mathbf{q})d\mathbf{q}\right)-\left(\dfrac{1}{2}\boldsymbol{\psi}(\mathbf{p})-\int_{\partial\Omega}\mathbf{n}_{\mathbf{p}}\cdot\mathbf{D}\nabla_{\mathbf{p}}G(\mathbf{q};\mathbf{p})\boldsymbol{\psi}(\mathbf{q})ds_{\mathbf{q}}\right)$ as the residual. \\
	
	Volume integral:\\
	$\bullet$ Get the volume integral $\mathbf{u}_{\text{h,V}}$ from solving the interface problem  (\ref{eq:volume}) by an FFT based Poisson solver.\\
	
	$\bullet$ Compute the limit values of the volume integral 
	\[ \dfrac{\partial\mathbf{u}_{\textrm{h,V}}^{+}}{\partial\mathbf{n}}=\mathbf{n_{p}}\cdot\mathbf{D}\nabla_{\mathbf{p}}\left(\int_{\Omega}G(\mathbf{q};\mathbf{p})\mathbf{f}(\mathbf{q})d\mathbf{q}\right)
	\]
	at the discretization nodes of the boundary curve by quadratic interpolation.  
	
	\While{$||\mathbf{r}||>tol_{\textrm{Richardson}}$}
	{
		
		Boundary integral:\\
		$\bullet$ Evaluate the jumps of the function and its partial derivatives 
		at the intersections to modify the right hand side of the discrete linear 
		system (\ref{eq:correct}) of the equivalent interface problem.
		
		$\bullet$ Get the boundary integral $\mathbf{u}_{\text{h,B}}$ from solving the interface problem (\ref{eq:single}) by an FFT based Poisson solver.
		
		$\bullet$ Compute the limit values of the boundary integral 
		\[ \dfrac{\partial\mathbf{u}_{\textrm{h,B}}^{+}}{\partial\mathbf{n}}=\dfrac{1}{2}\boldsymbol{\psi}(\mathbf{p})-\int_{\partial\Omega}\mathbf{n}_{\mathbf{p}}\cdot\mathbf{D}\nabla_{\mathbf{p}}G(\mathbf{q};\mathbf{p})\boldsymbol{\psi}(\mathbf{q})ds_{\mathbf{q}}
		\]
		at the discretization nodes of the boundary curve by quadratic interpolation.
		
		Compute the norm of the residual $\mathbf{r}$.
		
		Update the solution
		$\mathbf{u}_{\nu,\text{h}}^{n+1} \leftarrow \mathbf{u}_{\text{h,V}}+\mathbf{u}_{\text{h,B}}.
		$
		
		Update the iteration number $\nu \leftarrow \nu +1$.
		
	}
	
	{
		Accept the solution $\mathbf{u}_{\nu,\text{h}}^{n+1}$
		as an approximation of the vector function $\mathbf{u}$ at time $t=t^{n+1}$.
	}
	
\end{algorithm}

In the implementation, we first discretize the domain boundary $\partial\Omega$ 
by a set of quasi-uniformly spaced nodes and compute the normals, 
tangents, and curvatures of the curve or the surface at the 
discretization nodes of the domain boundary. We then partition 
the larger regular domain $\mathcal{B}$ into a uniform Cartesian 
grid and identify the regular and irregular grid nodes.
Locate all intersections of the interface with the Cartesian grid lines
and compute the unit normals and tangents at those intersections too.
All the above information is useful when computing the jumps 
of the function and its partial derivatives at the intersections 
and the discretization nodes of the domain boundary. Since the BIE 
(\ref{eq:bie}) to be solved can be noted as $\mathcal{A}\boldsymbol{\psi}=\mathbf{g}$, 
the simple Richardson iteration method shown in detailed algorithm 
\ref{alg:algo} is easily extended to the GMRES iteration.

\section{Numerical examples}

\begin{figure}[!t]
	\centering
	\includegraphics[scale=.5]{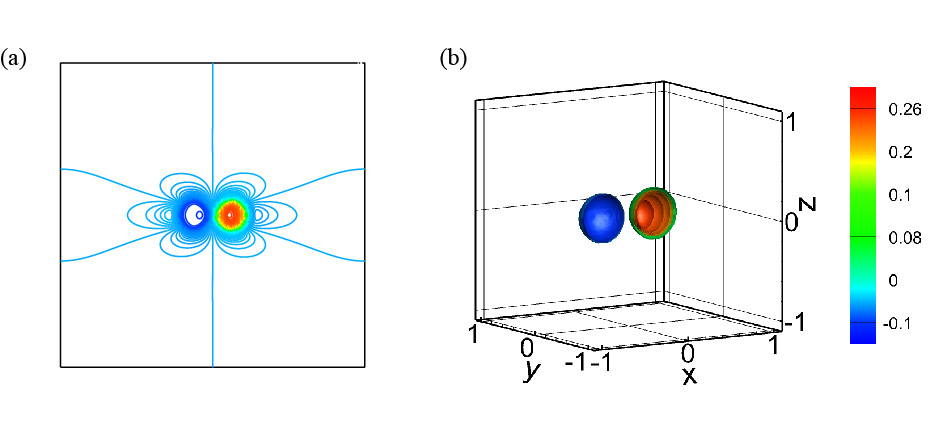}
	\caption{Extracellular stimulus expressed as a virtual
		electrode which consists of an anode and a cathode (a) in two-dimensional space and (b) in three-dimensional space, respectively.\label{fig:AC}}
\end{figure}

In this section, we present the numerical results from the KFBI method 
in different computational domains whose domain boundaries are 
represented by a set of $M_{\partial\Omega}$ quasi-uniform spaced 
nodes. Set the larger regular domain $\mathcal{B}$ as a square or 
a cube which is partitioned into a uniform Cartesian grid and 
denote $h$ as the mesh parameter. The transmembrane polarization 
is induced by the application of a virtual electrode consisting 
of a virtual cathode and a virtual anode \cite{vigmond2008solvers} 
and is given by

\[
I_{\text{stim}}(\mathbf{x})=\begin{cases}
V_{\text{stim}}^{\text{m}} & \textrm{if }\;r<0.25,\\
-V_{\text{stim}}^{\text{m}} & \textrm{if }\;s<0.25,
\end{cases}
\]
where $V_{\text{stim}}^{\text{m}}$ is a constant, $r$ is the 
distance between $\mathbf{x}$ and the anode, and $s$ is the
distance between $\mathbf{x}$ and the cathode. A depolarization 
after the application of the virtual battery in two- and three- 
dimensional spaces are shown in Fig. \ref{fig:AC}.

In this work, we use the simple but widely used FitzHugh-Nagumo 
(FHN) model as the membrane dynamics which has one gating variable 
and is governed by 

\[
I_{\text{ion}}(V_{\text{m}},q)=H\left(q-V_{\text{m}}(V_{\text{m}}-\theta)(1-V_{\text{m}})\right),
\]
and 

\[
\mathcal{M}(V_{\text{m}},q)=\alpha V_{\text{m}}-\zeta q,
\]
where $H=100$, $\theta=0.25$, $\alpha=5/20$, and $\zeta=1$. The gating 
variable $q$ and the transmembrane voltage $V_{\text{m}}$ are set to be at 
rest at time $t=0$. We model the larger regular domain $\mathcal{B}$ as 
consisting of fibers that lie parallel to the $x$-axis.

For all the numerical simulations, the membrane capacitance per unit area,
the surface-to-volume ratio, and the value of the extracellular stimulus
are set as $C_{\text{m}} = 1$, $\beta = 1000$, and $V_{\text{stim}}^{\text{m}}=10$,
respectively. The relative tolerances of the simple Richardson iteration and the 
GMRES iteration are chosen to be $10^{-8}$ and the iteration parameter $\gamma$
for the simple Richardson iteration is selected as 0.8. The positions of the
cathode and anode are set as $(\pm 0.3,0)$ in two-dimensional space and 
$(\pm 0.3, 0, 0)$ in three-dimensional space, respectively. The absolute 
tolerance in the Newton iteration method for solving the nonlinear systems
(\ref{eq:ode1}) and (\ref{eq:ode2}) is chosen to be $10^{-10}$. In all numerical
experiments, set the timestep size is the same as the mesh parameter, $i.e.$,
$\Delta t = h$. The simulations of the bidomain equations are implemented 
with custom codes written in C++ and the numerical simulations are all
performed on a 3.6 GHz computer with an Intel Core i3-4160 CPU.

In two-dimensional space, the bidomain equations are numerically simulated 
on a circular domain and a heart-shaped domain whose boundaries are presented
by a set of quasi-uniformly distributed nodes and the periodic cubic
splines on these points. In three-dimensional space, we only use the 
information of the discrete triangle grid of the domain boundary to 
simulate the bidomain equations by constructing local coordinate 
systems, instead of knowing the implicit expression of the domain boundary.
In this work, we solve the FHN bidomain equations in a spherical region.
The numerical results of the bidomain equations simulated on a human LV space
model are shown to demonstrate that the KFBI method 
works well for complex domains.

We first show the second-order convergence of the KFBI method.
In two-dimensional space, we take the solution of the bidomain equation
from the KFBI method simulated on a $2048\times2048$ Cartesian grid
with $M_{\partial\Omega} = 2048$ as the high precision solution. 
To numerically verify the second-order convergence of the KFBI method 
where a center finite difference method is applied, we estimate the 
following quality
\begin{equation}
\rho=\dfrac{log(e_{1}/e_{2})}{log(h_{1}/h_{2})},
\end{equation}
where $h_{1}$ and $h_{2}$ are two consecutive mesh parameters and $e_{1}$ 
and $e_{2}$ are the corresponding errors in the scaled $l^2$-norm or the 
infinity norm. Here, the scaled discrete $l^{2}-$norm of a vector
$\mathbf{v}=(v_{0},v_{1},\cdots,v_{n})^{T}\in\mathbb{R}^{n}$ is defined
by 
\begin{equation}
||\mathbf{v}||_{l^{2}}\equiv\sqrt{\dfrac{1}{n}\sum_{i=1}^{n}v_{i}^{2}}.
\end{equation}
As shown in Tables \ref{tab:error-fhn-circle-0.5}-\ref{tab:error-fhn-circle-2}
and Tables \ref{tab:error-fhn-heart-0.5}-\ref{tab:error-fhn-heart-2}, 
in two-dimensional space, the KFBI method has a second-order convergence 
obtained from the results solved on a circular domain and a heart-shaped 
domain. In three-dimensional space, we take the solution of the bidomain 
equation from the KFBI method simulated on a $512\times512\times512$ 
Cartesian grid with $M_{\partial\Omega} = 2052$ as the high precision 
solution. As shown in in Tables \ref{tab:error-fhn-3d-0.25}-\ref{tab:error-fhn-3d-0.75},
the numerical results demonstrate that the KFBI can also achieve a 
second-order convergence in three-dimensional space, which agree 
with our expectations well. In Tables \ref{tab:error-fhn-circle-2}
and \ref{tab:error-fhn-heart-2}, the last two columns show the average
iteration numbers of the GMRES iteration in each timestep and
the total CPU times (in seconds) at time $t=2$.



We then give the membrane voltage $V_{\text{m}}$ simulated from the KFBI method  
on a $256\times256$ and $128\times128\times128$ Cartesian grid, 
respectively. Figs. \ref{fig:fhn-circle-isoplots} and 
\ref{fig:fhn-heart-isoplots} 
show the contours for the transmembrane voltage $V_{\text{m}}$ 
on different computational domains. The transmembrane voltage at 
discrete times after the application of the virtual battery 
in a sphere domain are shown in Fig. \ref{fig:sphere_t=00003D1.5-2.5}. 
The trajectories of the transmembrane voltage collected at different 
points and the iteration number of the GMRES iteration or the simple 
Richardson iteration are shown.

We finally show the numerical results simulated on a real human 
LV model which was constructed from an vivo magnetic resonance
imaging study of a healthy volunteer \cite{Cai2015LV}. 
As shown in Fig. , we only need the information of the discretization
nodes in the discrete triangle grid of the domain boundary.
As shown in Fig. \ref{fig:lv}, we stimulate one corner of the LV 
model at the very beginning without applying any extracellular 
stimulus and the numerical results demonstrate that the KFBI method works 
in very complex domains as well as in the easy cases.
Scroll wave results of the FHN bidomain equation in the LV space model 
are shown in Fig. \ref{fig:scroll-lv}.

\subsection{Example 1}
\begin{figure}[!t]
	\centering
	\includegraphics[scale=.33]{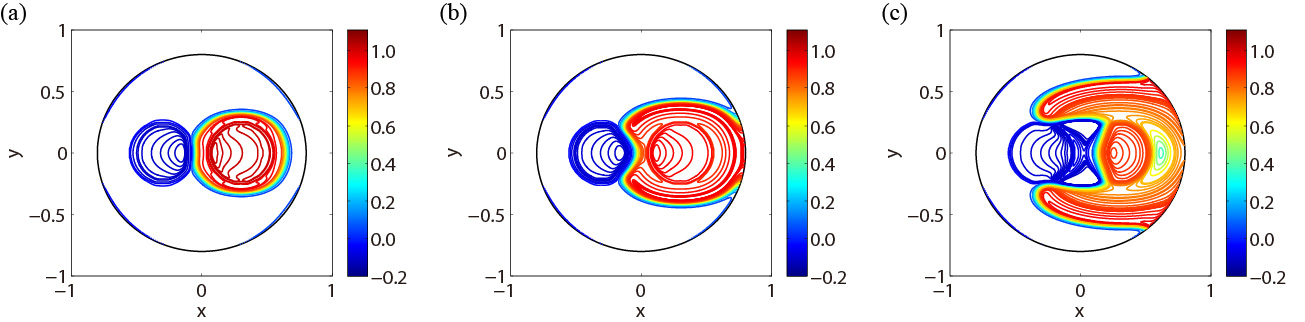}
	\caption{Iso-contours of the transmembrane voltage solved 
		on the circular domain with a 256$\times$256 Cartesian grid at time $t=0.5$,
		$t=1$, and $t=2$.\label{fig:fhn-circle-isoplots}}
\end{figure}

\begin{figure}[!t]
	\centering
	\includegraphics[scale=.5]{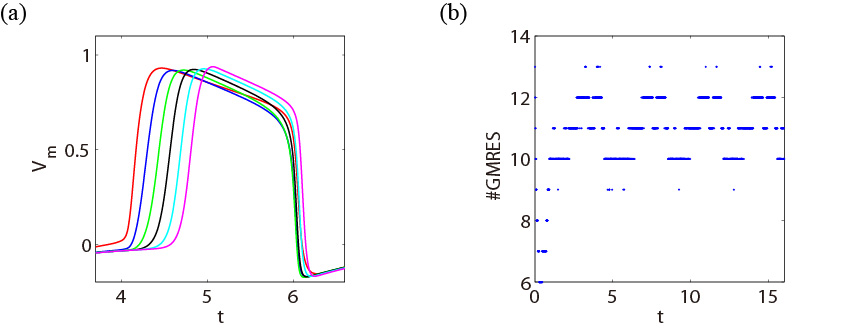}
	\caption{Results of FHN bidomain model on the circular domain with a 256$\times$256
		Cartesian grid. (a) Trajectories of the transmembrane voltage collected at six points distributed
		on a straight line from $(0.5625,0)$ to $(0.71875,0)$ evenly. 
		(b) The iteration number of GMRES iteration against time. \label{fig:fhn-circle-vt}}
\end{figure}

In the case of the two-dimensional space, the intracellular 
conductivities are 30 along the fiber direction and 5 perpendicular to 
the fiber direction, while the extracellular conductivities are 20 along 
the fiber direction and 10 perpendicular to the fiber direction, $i.e.$, 
$\sigma_{\text{i}}=30$, $\eta_{\text{i}}=5$,
$\sigma_{\text{e}}=20$ and $\eta_{\text{e}}=10$.
Tables \ref{tab:error-fhn-circle-0.5}-\ref{tab:error-fhn-circle-2}
and Figs. \ref{fig:fhn-circle-isoplots} and \ref{fig:fhn-circle-vt}
show the numerical results  simulated on a circular domain, $e.g.$, 
the convergence order, CPU cost times (seconds), iso-contours of 
the transmembrane voltage, and the iteration numbers of the
GMRES iteration.

\begin{table*}[!t]
	\centering
	\begin{tabular}{|c|c|c|c|c|c|}
		\hline 
		$M_{\partial\Omega}$ & Grid Size & $||e_{h}||_{l^{2}}$ & $\rho$ & $||e_{h}||_{\infty}$ & $\rho$\tabularnewline
		\hline 
		\hline 
		64 & 64$\times$64 & 0.063262 & - & 0.663016 & -\tabularnewline
		\hline 
		128 & 128$\times$128 & 0.016879 & 1.91 & 0.209426 & 1.66\tabularnewline
		\hline 
		256 & 256$\times$256 & 0.004422 & 1.93 & 0.061362 & 1.77\tabularnewline
		\hline 
		512 & 512$\times$512 & 0.001150 & 1.94 & 0.015737 & 1.96\tabularnewline
		\hline 
		1024 & 1024$\times$1024 & 0.000313 & 1.88 & 0.004245 & 1.89\tabularnewline
		\hline 
	\end{tabular}
	\caption{Errors of the transmembrane voltage in the scaled discrete
		$l^2$-norm and infinity norm solved on the circular domain at time t=0.5.
		\label{tab:error-fhn-circle-0.5}}
\end{table*}

\begin{table*}[!t]
	\centering
	\begin{tabular}{|c|c|c|c|c|c|}
		\hline 
		$M_{\partial\Omega}$ & Grid Size & $||e_{h}||_{l^{2}}$ & $\rho$ & $||e_{h}||_{\infty}$ & $\rho$\tabularnewline
		\hline 
		\hline 
		64 & 64$\times$64 & 0.156928 & - & 0.892311 & -\tabularnewline
		\hline 
		128 & 128$\times$128 & 0.037347 & 2.07 & 0.371083 & -\tabularnewline
		\hline 
		256 & 256$\times$256 & 0.008594 & 2.12 & 0.092505 & 2.00\tabularnewline
		\hline 
		512 & 512$\times$512 & 0.002115 & 2.02 & 0.022442 & 2.04\tabularnewline
		\hline 
		1024 & 1024$\times$1024 & 0.000510 & 2.05 & 0.005078 & 2.14\tabularnewline
		\hline 
	\end{tabular}
	\caption{Errors of the transmembrane voltage in the scaled discrete
		$l^2$-norm and infinity norm solved on the circular domain at time t=1.\label{tab:error-fhn-circle-1}}
\end{table*}

\begin{table*}[!t]
	\centering
	\begin{tabular}{|c|c|c|c|c|c|}
		\hline 
		$M_{\partial\Omega}$ & Grid Size & $||e_{h}||_{l^{2}}$ & $\rho$ & $||e_{h}||_{\infty}$ & $\rho$\tabularnewline
		\hline 
		\hline 
		64 & 64$\times$64 & 0.243811 & - & 0.943763 & -\tabularnewline
		\hline 
		128 & 128$\times$128 & 0.062995 & 1.95 & 0.550973 & -\tabularnewline
		\hline 
		256 & 256$\times$256 & 0.013931 & 2.18 & 0.127781 & 2.11\tabularnewline
		\hline 
		512 & 512$\times$512 & 0.003317 & 2.07 & 0.030559 & 2.06\tabularnewline
		\hline 
		1024 & 1024$\times$1024 & 0.000734 & 2.18 & 0.006604 & 2.21\tabularnewline
		\hline 
	\end{tabular}
	\caption{Errors of the transmembrane voltage in the scaled discrete
		$l^2$-norm and infinity norm solved on the circular domain at time t=1.5.\label{tab:error-fhn-circle-1.5}}
\end{table*}

\begin{table*}[!t]
	\centering
	\begin{tabular}{|c|c|c|c|c|c|c|c|}
		\hline 
		$M_{\partial\Omega}$ & Grid Size & $||e_{h}||_{l^{2}}$ & $\rho$ & $||e_{h}||_{\infty}$ & $\rho$ & \#GMRES & CPU (s)\tabularnewline
		\hline 
		\hline 
		64 & 64$\times$64 & 0.312412 & - & 0.954275 & - & 10.03 & 4.98\tabularnewline
		\hline 
		128 & 128$\times$128 & 0.085055 & 1.88 & 0.698887 & - & 9.65 & 35.70\tabularnewline
		\hline 
		256 & 256$\times$256 & 0.018857 & 2.17 & 0.184695 & 1.92 & 8.94 & 272.12\tabularnewline
		\hline 
		512 & 512$\times$512 & 0.004393 & 2.10 & 0.043494 & 2.09 & 8.07 & 2127.33\tabularnewline
		\hline 
		1024 & 1024$\times$1024 & 0.000948 & 2.21 & 0.009010 & 2.27 & 7.27 & 16696.65\tabularnewline
		\hline 
	\end{tabular}
	\caption{Errors of the transmembrane voltage in the scaled discrete
		$l^2$-norm and infinity norm solved on the circular domain at time t=2.\label{tab:error-fhn-circle-2}}
\end{table*}

\subsection{Example 2}
\begin{figure}[!t]
	\centering
	\includegraphics[width=1\textwidth]{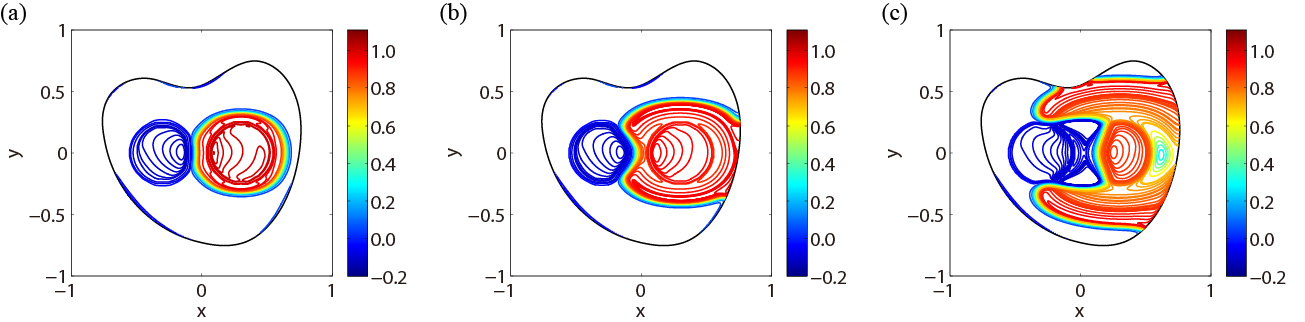}
	\caption{Iso-contours of the transmembrane voltage solved
		on the heart-shaped domain with a 256$\times$256 Cartesian grid at time $t=0.5$,
		$t=1$ and $t=2$.\label{fig:fhn-heart-isoplots}}
\end{figure}

\begin{figure}[!t]
	\centering
	\includegraphics[width=1\textwidth]{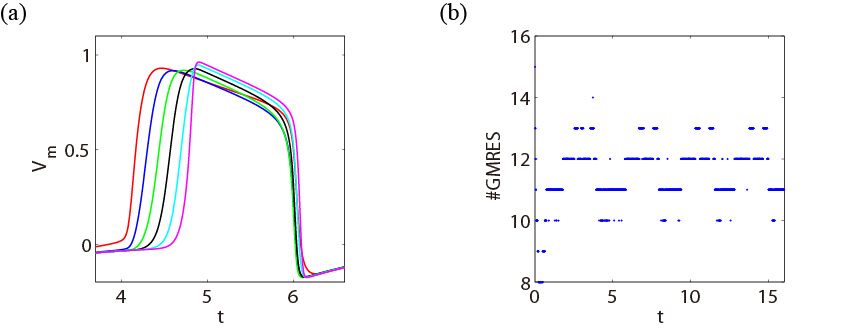}
	\caption{Results of FHN bidomain model on the heart-shaped domain with a 256$\times$256
		Cartesian grid. (a) Trajectories of the transmembrane voltage collected at six points distributed
		on a straight line from $(0.5625,0)$ to $(0.71875,0)$ evenly. (b) The iteration number of GMRES iteration against time. \label{fig:fhn-heart-vt}}	
\end{figure}

Simulating the bidomain equation on the heart-shaped domain
in the two-dimensional space, the intracellular and extracellular 
conductivities are the same as noted previously.
Tables \ref{tab:error-fhn-heart-0.5}-\ref{tab:error-fhn-heart-2}
and Figs. \ref{fig:fhn-heart-isoplots} and \ref{fig:fhn-heart-vt}
show the numerical results  simulated on a circular domain, $e.g.$, 
the convergence order, CPU cost times (seconds), iso-contours of 
the transmembrane voltage, and the iteration numbers of the
GMRES iteration.

\begin{table*}[!t]
	\centering
	\begin{tabular}{|c|c|c|c|c|c|}
		\hline 
		$M_{\partial\Omega}$ & Grid Size & $||e_{h}||_{l^{2}}$ & $\rho$ & $||e_{h}||_{\infty}$ & $\rho$\tabularnewline
		\hline 
		\hline 
		64 & 64$\times$64 & 0.068323 & - & 0.669969 & -\tabularnewline
		\hline 
		128 & 128$\times$128 & 0.018314 & 1.90 & 0.227030 & 1.56\tabularnewline
		\hline 
		256 & 256$\times$256 & 0.004805 & 1.93 & 0.062225 & 1.87\tabularnewline
		\hline 
		512 & 512$\times$512 & 0.001253 & 1.94 & 0.015830 & 1.97\tabularnewline
		\hline 
		1024 & 1024$\times$1024 & 0.000340 & 1.88 & 0.004270 & 1.89\tabularnewline
		\hline 
	\end{tabular}
	\caption{Errors of the transmembrane voltage in the scaled discrete
		$l^2$-norm and infinity norm solved on the heart-shaped domain at time t=0.5.\label{tab:error-fhn-heart-0.5}}	
\end{table*}

\begin{table*}[!t]
	\centering
	\begin{tabular}{|c|c|c|c|c|c|}
		\hline 
		$M_{\partial\Omega}$ & Grid Size & $||e_{h}||_{l^{2}}$ & $\rho$ & $||e_{h}||_{\infty}$ & $\rho$\tabularnewline
		\hline 
		\hline 
		64 & 64$\times$64 & 0.169297 & - & 1.037716 & -\tabularnewline
		\hline 
		128 & 128$\times$128 & 0.040869 & 2.05 & 0.381248 & 1.44\tabularnewline
		\hline 
		256 & 256$\times$256 & 0.009457 & 2.11 & 0.130960 & 1.54\tabularnewline
		\hline 
		512 & 512$\times$512 & 0.002317 & 2.03 & 0.023010 & 2.51\tabularnewline
		\hline 
		1024 & 1024$\times$1024 & 0.000544 & 2.09 & 0.005196 & 2.15\tabularnewline
		\hline 
	\end{tabular}
	\caption{Errors of the transmembrane voltage in the scaled discrete
		$l^2$-norm and infinity norm solved on the heart-shaped domain at time t=1.\label{tab:error-fhn-heart-1}}
\end{table*}

\begin{table*}[!t]
	\centering
	\begin{tabular}{|c|c|c|c|c|c|}
		\hline 
		$M_{\partial\Omega}$ & Grid Size & $||e_{h}||_{l^{2}}$ & $\rho$ & $||e_{h}||_{\infty}$ & $\rho$\tabularnewline
		\hline 
		\hline 
		64 & 64$\times$64 & 0.261088 & - & 0.950014 & -\tabularnewline
		\hline 
		128 & 128$\times$128 & 0.068827 & 1.92 & 0.546080 & 0.80\tabularnewline
		\hline 
		256 & 256$\times$256 & 0.015225 & 2.18 & 0.130868 & 2.06\tabularnewline
		\hline 
		512 & 512$\times$512 & 0.003636 & 2.07 & 0.031317 & 2.06\tabularnewline
		\hline 
		1024 & 1024$\times$1024 & 0.000803 & 2.18 & 0.006758 & 2.21\tabularnewline
		\hline 
	\end{tabular}
	\caption{Errors of the transmembrane voltage in the scaled discrete
		$l^2$-norm and infinity norm solved on the heart-shaped domain at time t=1.5.\label{tab:error-fhn-heart-1.5}}
\end{table*}

\begin{table*}[!t]
	\centering
	\begin{tabular}{|c|c|c|c|c|c|c|c|}
		\hline 
		$M_{\partial\Omega}$ & Grid Size & $||e_{h}||_{l^{2}}$ & $\rho$ & $||e_{h}||_{\infty}$ & $\rho$ & \#GMRES & CPU (s)\tabularnewline
		\hline 
		\hline 
		64 & 64$\times$64 & 0.336369 & - & 0.970956 & - & 14.23 & 6.81\tabularnewline
		\hline 
		128 & 128$\times$128 & 0.090293 & 1.90 & 0.754806 & - & 11.36 & 41.39\tabularnewline
		\hline 
		256 & 256$\times$256 & 0.020045 & 2.17 & 0.214061 & 1.82 & 12.49 & 354.83\tabularnewline
		\hline 
		512 & 512$\times$512 & 0.004666 & 2.10 & 0.053079 & 2.01 & 9.54 & 2388.28\tabularnewline
		\hline 
		1024 & 1024$\times$1024 & 0.001004 & 2.22 & 0.011413 & 2.22 & 8.77 & 19094.92\tabularnewline
		\hline 
	\end{tabular}
	\caption{Errors of the transmembrane voltage in the scaled discrete
		$l^2$-norm and infinity norm solved on the heart-shaped domain at time t=2.\label{tab:error-fhn-heart-2}}
\end{table*}

\subsection{Example 3}
\begin{figure}[!t]
	\centering
	\includegraphics[width=1\textwidth]{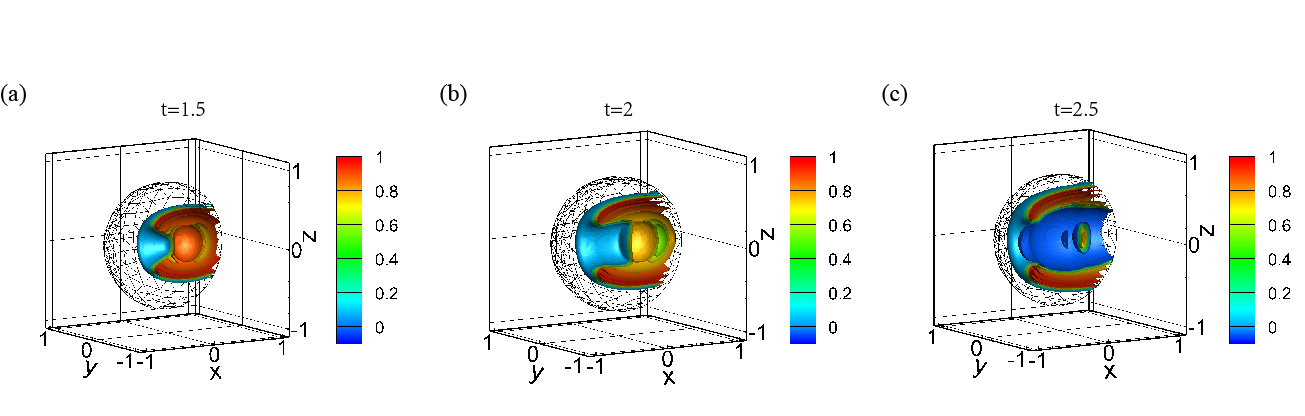}
	\centering{}\caption{Iso-contours of the transmembrane voltage solved 
		in the spherical region with a 128$\times$128$\times$128
		Cartesian grid at discrete time $t=1.5$, $t=2$, and $t=2.5$. 
		The observing domain is limited with
		$y\geq0$.\label{fig:sphere_t=00003D1.5-2.5}}
\end{figure}

\begin{figure}[!t]
	\centering
	\includegraphics[width=1\textwidth]{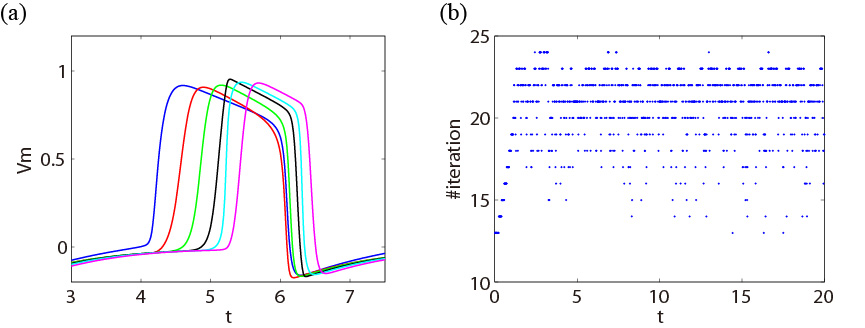}
	\caption{Results of FHN bidomain model in the spherical region with a 128$\times$128$\times$128
		Cartesian grid. (a) Trajectories of the transmembrane voltage collected at six points distributed
		on a straight line from $(0.5625,0,0)$ to $(0.875,0,0)$ evenly. (b) The number of the simple 
		Richardson iteration against
		time and the average iteration number is 21.\label{fig:3d-t-itr}}	
\end{figure}

In this subsection, we show the numerical results of the bidomain equation 
simulated in a sphere domain. In the case of the three-dimensional case, 
set the conductivity in the direction orthogonal to the fiber axis be the same 
as the conductivity in the direction perpendicular to the fiber, $i.e.$, 
$\mathbf{D}_{\text{i}}=\text{diag}(30,5,5)$ and $\mathbf{D}_{\text{i}}=\text{diag}(20,10,10)$.
Tables \ref{tab:error-fhn-3d-0.25}-\ref{tab:error-fhn-3d-0.75}
and Figs. \ref{fig:sphere_t=00003D1.5-2.5} and \ref{fig:3d-t-itr}
show the numerical results  simulated on a circular domain, $e.g.$, 
the convergence order, CPU cost times (seconds), iso-contours of 
the transmembrane voltage, and the iteration numbers of the
simple Richardson iteration.


\begin{table*}[!t]
	\centering
	\begin{tabular}{|c|c|c|c|c|c|}
		\hline 
		$M_{\partial\Omega}$ & Grid Size & $||e_{h}||_{l^{2}}$ & $\rho$ & $||e_{h}||_{\infty}$ & $\rho$\tabularnewline
		\hline 
		\hline 
		18 & 32$\times$32$\times$32 & 4.16e-02 & - & 5.25e-01 & -\tabularnewline
		\hline 
		66 & 64$\times$64$\times$64 & 1.34e-02 & 1.6 & 2.12e-01 & 1.3\tabularnewline
		\hline 
		258 & 128$\times$128$\times$128 & 3.04e-03 & 2.1 & 7.06e-02 & 1.6\tabularnewline
		\hline 
		1026 & 256$\times$256$\times$256 & 6.20e-04 & 2.3 & 2.29e-02 & 1.6\tabularnewline
		\hline 
	\end{tabular}
	\caption{Errors of the transmembrane voltage in the scaled discrete
		$l^2$-norm and infinity norm solved in the spherical region at time t=0.25.\label{tab:error-fhn-3d-0.25}}
\end{table*}

\begin{table*}[!t]
	\centering
	\begin{tabular}{|c|c|c|c|c|c|}
		\hline 
		$M_{\partial\Omega}$ & Grid Size & $||e_{h}||_{l^{2}}$ & $\rho$ & $||e_{h}||_{\infty}$ & $\rho$\tabularnewline
		\hline 
		\hline 
		18 & 32$\times$32$\times$32 & 4.52e-02 & - & 8.00e-01 & -\tabularnewline
		\hline 
		66 & 64$\times$64$\times$64 & 1.82e-02 & 1.3 & 4.92e-01 & 0.7\tabularnewline
		\hline 
		258 & 128$\times$128$\times$128 & 3.16e-03 & 2.5 & 1.37e-01 & 1.8\tabularnewline
		\hline 
		1026 & 256$\times$256$\times$256 & 6.11e-04 & 2.4 & 2.79e-02 & 2.3\tabularnewline
		\hline 
	\end{tabular}
	\caption{Errors of the transmembrane voltage in the scaled discrete
		$l^2$-norm and infinity norm solved
		in the spherical region at time t=0.5.\label{tab:error-fhn-3d-0.5}}
\end{table*}

\begin{table*}[!t]
	\centering
	\begin{tabular}{|c|c|c|c|c|c|c|c|}
		\hline 
		$M_{\partial\Omega}$ & Grid Size & $||e_{h}||_{l^{2}}$ & $\rho$ & $||e_{h}||_{\infty}$ & $\rho$ & \#SIMPLE & CPU (s)\tabularnewline
		\hline 
		\hline 
		18 & 32$\times$32$\times$32 & 8.69e-02 & - & 9.11e-01 & - & 18 & 22.4\tabularnewline
		\hline 
		66 & 64$\times$64$\times$64 & 2.25e-02 & 1.9 & 5.13e-01 & 0.8 & 19 & 195.2\tabularnewline
		\hline 
		258 & 128$\times$128$\times$128 & 4.30e-03 & 2.4 & 1.20e-01 & 2.1 & 21 & 3296.0\tabularnewline
		\hline 
		1026 & 256$\times$256$\times$256 & 7.27e-04 & 2.6 & 2.49e-02 & 2.3 & 16 & 62291.2\tabularnewline
		\hline 
	\end{tabular}
	\caption{Errors of the transmembrane voltage in the scaled discrete
		$l^2$-norm and infinity norm solved
		in the spherical region at time t=0.75.\label{tab:error-fhn-3d-0.75}}
\end{table*}

\subsection{Example 4}

\begin{figure}[!t]
	\centering
	\includegraphics[scale=.5]{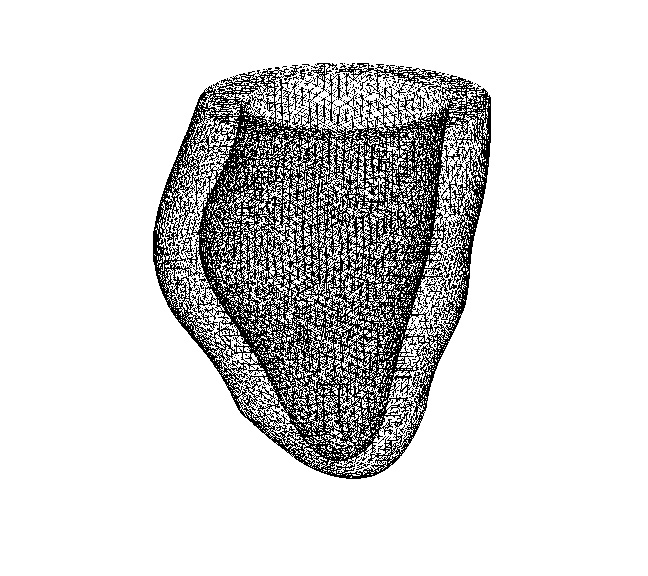}
	\caption{Illustration of the discrete triangle grid of boundary of the LV spatial model.\label{fig:lv_model}}	
\end{figure}

\begin{figure}[!t]
	\centering
	\includegraphics[width=1\textwidth]{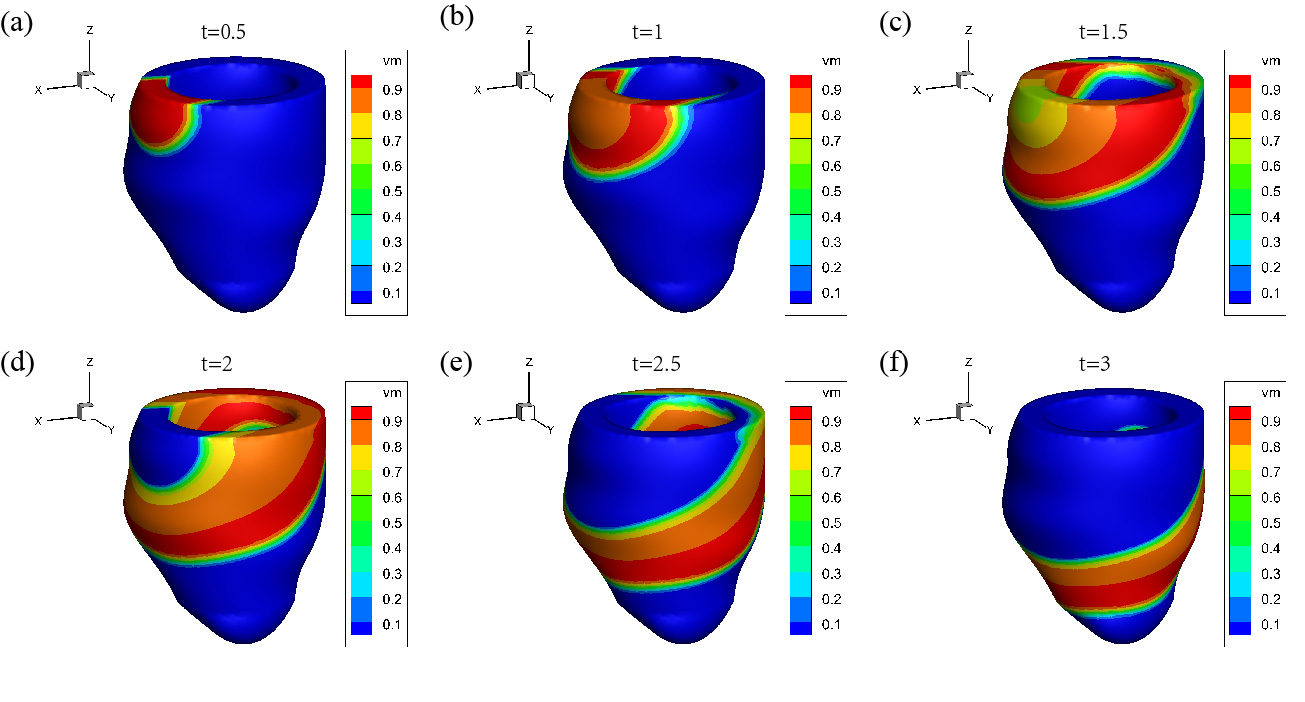}
	\caption{Iso-contours of the transmembrane voltage solved
		on the human LV model with a 128$\times$128$\times$128 Cartesian grid at discrete
		times.\label{fig:lv}}	
\end{figure}

\begin{figure}[!t]
	\centering
	\includegraphics[width=1\textwidth]{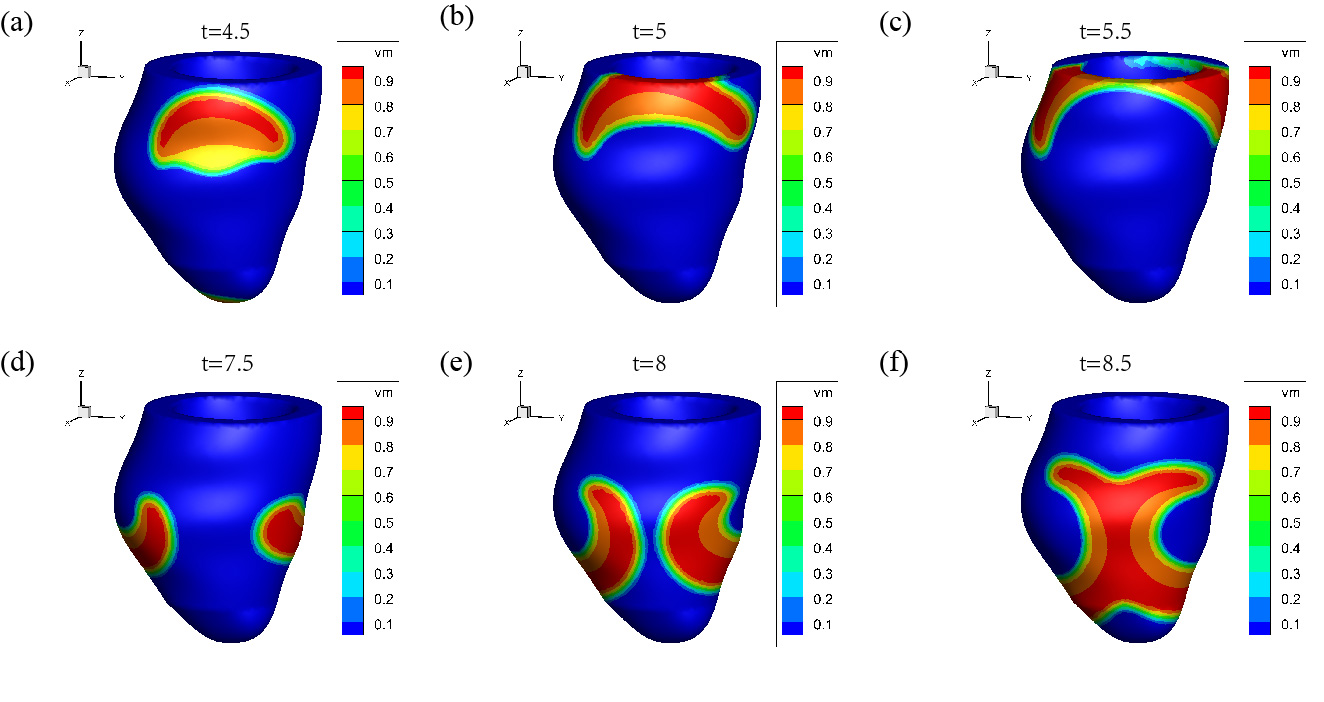}
	\caption{scroll wave results of the FHN bidomain equations simulated on the human LV
		model.\label{fig:scroll-lv}}
\end{figure}

We finally simulate the FHN bidomain equations on a LV model
which was constructed from an vivo magnetic resonance
imaging study of a healthy volunteer \cite{Cai2015LV}. 
As shown in Fig. \ref{fig:lv_model}, in three-dimensional space, we construct
local coordinate systems on the discrete triangle grid of
the domain boundary to calculate the normals, tangents, 
and curvatures of the surface at the discretization nodes
of the domain boundary. The initial conditions of the 
transmembrane voltage $V_{\text{m}}$ and the stage variable 
$q$ are given by the smooth data
\[
V_{\text{m}}(x,y,z,0)=\begin{cases}
1 & \textrm{if }\sqrt{(x-0.3)^{2}+(y-0.05)^{2}+(z-0.35)^{2}}<0.1,\\
0 & \textrm{otherwise},
\end{cases}
\]
and 

\[
q(x,y,z,0)=0.
\]
Set the larger regular domain $\mathcal{B}$ to be a cube $[-0.4,0.4]^{3}$. 
The parameters, tolerances, and modules are set the same as those noted 
previously in this work except that the iteration parameter $\gamma$ for 
the simple Richardson iteration method is 0.7.
A second stimulus 

\[
I_{s}(x,y,z)=\begin{cases}
1 & \textrm{if }\sqrt{(x-0.25)^{2}+(y-0.05)^{2}+(z-0.05)^{2}}<0.15,\\
0 & \textrm{otherwise}
\end{cases}
\]
is applied to the transmembrane voltage at time $t=3.5$. Figs. \ref{fig:lv} 
and \ref{fig:scroll-lv} show the traveling waves and scroll
wave results simulated on the human LV model.

\section{Discussion}

The bidomain equations which mathematically model the electrical activity
of the cardiac tissue are space- and time-dependent problems. In this work, we
use the second-order Strang splitting to separate the linear diffusion
part from the nonlinear reaction part and integrate the corresponding 
diffusion equation by the implicit midpoint scheme. The resulting 
semi-discrete linear diffusion equation which is a homogeneous Neumann 
boundary value problem can be efficiently solved by the KFBI method 
without the requirement of the analytical expression of Green's functions.
In this work, the KFBI method uses the second-order centered difference 
method to spatially discrete the equivalent interface problems and uses 
the quadratic interpolation to calculate the limit values of the integrals 
at the discretization nodes of the domain boundary. At the same time, the 
second-order expressions of the domain boundary and a second-order Taylor 
expansion to compute the jumps of the function and its derivatives across
the domain boundary are used in the KFBI method. Numerical results simulated on 
different domains demonstrate that the KFBI method has second-order accuracy
and the iteration number of the simple Richardson iteration and the GMRES 
iteration is independent of the mesh parameter of the Cartesian grid. 

The KFBI method is suitable for the complex computational domains. After 
embedding the complex computational domain into a larger regular domain, 
according to the potential theory, the boundary value problems obtained 
from the linear diffusion part in the complex computational domain can 
be reformulated into the second kind of the Fredholm BIEs on the larger 
regular domain. During each iteration for solving the BIE, the process
of solving the equivalent interface problems on the Cartesian grid of the
larger regular domain to calculate the corresponding integrals is easy to 
implement. Moreover, in three-dimensional space, we only need the discrete triangle
grid of the domain boundary to estimate the information of the boundary 
curve by constructing local coordinate systems. The numerical results 
demonstrate that the KFBI method works well simulating the bidomain equations
in a human LV spatial model. Another advantage of the KFBI method is that the 
coefficient matrix of the discrete linear system obtained from discretizing
the equivalent interface problem is still symmetric and positive definite after 
appropriate modification at the irregular grid nodes. Therefore, the
fast elliptic solvers such as an FFT based Poisson solver can be used to 
efficiently solve the discrete linear system.

We point out that the KFBI method is not limited to the second-order
convergence version which is co-determined by the centered difference 
method, the quadratic interpolation, the second-order expression of 
the domain boundary, and the second-order Taylor expansion.
The KFBI method can be readily extended to a higher-order version 
\cite{xie2019fourth,ying2018cartesian} as long as the discrete method
and the Taylor expansion involved in the KFBI method are raised to the 
corresponding order which involves more grid nodes and more Taylor 
expansion terms.  

We emphasis that not only the finite difference method but also the 
finite element method can be applied for the spatial discretization 
of the equivalent interface problems. The larger regular domain 
$\mathcal{B}$ can be chosen flexibly to be a triangle, a rectangle, 
a circle, or any other regularly shaped domains as long as the Green's 
function on $\mathcal{B}$ exists and the structured grid-based fast 
elliptic solvers are applicable. Based on the flexibility of the KFBI 
method, the bidomain equations can be simulated in more general cases, 
for example, the conductivity is anisotropic and spatially dependent
or the computational domain changes over time.



%
%
%


%
%
%
%
%
%
%
%

\bibliographystyle{unsrt}
\bibliography{reference}


\end{document}